\documentclass[a4paper,12pt]{article}
\usepackage{latexsym}
\usepackage{amsfonts}
\usepackage{amsmath}
\usepackage{amssymb}
\def \R {{\rm I\!R}}

\def \N {{\rm I\!N}}

\def \L {{\rm I\!L}}

\def \H {{\rm I\! H}}

\def \E {{\rm E}}

\title{Random fields, large deviations and triviality in quantum field theory. Part I.} 

\author{ Adnan Aboulalaa\footnotemark[1]}
\date{}

\begin{document}
\vskip2cm
\newtheorem{thm}{Theorem}[section]
\newtheorem{defi}[thm]{Definition}
\newtheorem{prop}[thm]{Proposition}
\newtheorem{Rk}[thm]{Remark}
\newtheorem{lm}{Lemma}[section]
\newtheorem{co}[thm]{Corollary}
\maketitle
\footnotetext[1]{  E-mail: adnan.aboulalaa@polytechnique.org }
\begin{abstract}
\noindent
The issue of the existence and possible triviality of the Euclidean quantum scalar field in dimension 4 is investigated by using some large deviations techniques. As usual, the field $\varphi_{d}^{4}$ is obtained as a limit of regularized fields $\varphi_{k}^{4}$ associated with
a probability measures $\mu_{k,V}$, where $k, V$ represent ultraviolet and volume cutoffs. The result obtained is that in a fixed volume, the almost sure limit (as $k \rightarrow \infty$) of the density of $\mu_{k,V}$, with respect to the Gaussian free field measure, exists and is equal to $0$, when the coupling constant is not vanishing. This implies that $\mu_{k,V}$ can not have a strong limit as the ultraviolet cutoff is removed. 
Furthermore, the normalization sequence $Z_{k,V}=\E e^{-{\cal A}_{k,V}}$ is divergent as $k \rightarrow \infty$ for dimensions $d\geq4$ when the vacuum renormalization is lower than some threshold, which leads to the non ultraviolet stability of the field in this case.
These assertions are also valid for vector fields and can be extended to polynomial Lagrangians.
\\ \\


\end{abstract}


\section{Introduction}

Let $d>0$ and ${\cal P}={\cal S}'(\R^{d})$ be the space of tempered distributions, and $\mu_{0}$ be the Gaussian measure on ${\cal P}$
associated to the free field. In this paper we are interested in the class of random fields corresponding to the probability measure $\mu$ 
given by:
\begin{equation}
\label{Intro1}
\frac{d\mu}{d\mu_{0}}= \frac{1}{Z} \exp (-\int_{\R^{d}}L(\varphi(x), \partial \varphi(x))dx), \; \;  Z =E_{\mu_{0}}e^{-\int_{\R^{d}} L(\varphi(x), \partial \varphi(x)) dx}  
\end{equation}
\noindent
where $L$ is a Lagrangian density, which, in general, is a polynomial function. The fact that the sample fields $\varphi$ are irregular distributions, 
so that neither the pointwise values nor their products are well defined, makes the expression (\ref{Intro1}) a formal one.
\\
The connection between (\ref{Intro1}) and field theory is not obvious at all and had taken quite a long time to be established. For a survey of the subject 
we refer to Jaffe \cite{Jaffe2007} and Summers \cite{Summers}. 
In very few words, let us say that Quantum Field Theory (QFT) is currently the theoretical framework of modern particle physics with predictions that agree
with experimental results with very high precision orders. 
That the mathematical foundation of this theory is problematic was recognized in the 1950s and has given rise to a new discipline in mathematical physics, 
where two approaches have emerged: the algebraic quantum field theory (Haag, Kastler, Araki, see, e.g., \cite{Haag}, \cite{Araki}) and what can be called the analytic 
approach (Wightman, G\aa rding, see \cite{Jost}, \cite{SW}, \cite{Str}, \cite{Klauder}, \cite{Dimock} for an introduction to these topics), 
in which we find most of the Constructive Quantum Field Theory program, whose aim is the construction of rigorous mathematical models of QFT. 
\\
\\
\noindent
The objects dealt with in QFT are operators $\hat{\varphi}(x)$ indexed by $\R^{4}$ which is considered as a Minkowski space ; these operators act on some Hilbert space $\H$. 
It was found that the operators $\hat{\varphi}(x)$ are not regular with respect 
to $x$ and a proposal has been made to consider them as operator valued distributions $\hat{\varphi}(f)$ indexed by a space of test functions $f$ and satisfying certain conditions, the G\aa rding-Wightman axioms. Furthermore as QFT considers also objects like $W(x_{1}, .... x_{p})= \langle \hat{\varphi}(x_{1})...\hat{\varphi}(x_{p})\Psi, \Psi \rangle_{\H}$, where $\Psi$ is a special state (the vaccum), it was proved that the operators $\hat{\varphi}(f)$ can be reconstructed from a given distributions $W(f_{1}, .... f_{p})$ (the Wightman functions) that satisfy a set of axioms.
\\
While the {\it axioms} are justified and accepted, the pivotal issue was (and still is) whether there exists a non trivial field that corresponds to 
those of theoretical physics. This matter has been extensively studied, cf. Glimm, Jaffe \cite{GJ87} for references ; see also the survey and
bibliography in Summers \cite{Summers} and Malyshev \cite{Malyshev} for an overview of the probabilistic aspects.
\\
\noindent
Early constructions have been performed in the operator framework in the 1960s. To begin with dimension 2, following a pioneering work by Nelson \cite{Nelson1}, a positive
answer was given for the scalar field by Glimm and Jaffe \cite{GJ68} in finite volume and subsequently in infinite volume (Glimm, Jaffe, Spencer \cite{GJS74}).
In the meantime, the Euclidean treatment of these problems was developed by many authors (Nelson, Symanzik \cite{Symanzik}, \cite{Nelson2}, \cite{Nelson3} ) and proved to be much more 
convenient than the Minkowski setting. Osterwalder and Schrader \cite{OS} discovered that the Wightman distributions $W$ can be associated to their Euclidean counterpart $S(f_{1}, .... f_{p})$, the Schwinger functions that fulfill a set of conditions, the Osterwalder-Schrader (OS) axioms, see Zinoviev \cite{Zinoviev} for the question of the equivalence of the two constructions. 
Even more, it was noticed (see, e.g., \cite{Fr74}) that these Euclidean fields can be constructed through a probability measure on the some space, which is often ${\cal S}'(\R^{d})$, and the $S(f_{1}, .... f_{p})$ are the moments of this probability measure ; we end at this point this very short and incomplete description of the link between QFT and probability measures like (\ref{Intro1}).
\\
After the construction of quantum fields in dimension 2, the case of dimension 3 was also solved by Glimm-Jaffe \cite{GJ73}, Feldman-Osterwalder \cite{FO}, and 
Magnen-S\'eneor \cite{MS} and many authors, with different methods (\cite{Park}, \cite{Benfatto-al}, \cite{Balaban}, \cite{BF}, \cite{BFS}, etc.).
\\
In all these constructions, the interacting field is obtained as a limit of regularized fields $\varphi_{k,V}$ with volume and cutoffs denoted by $V$ and $k$ ; for the later, 
a momentum or lattice regularizations are used. In the case of a momentum cutoff, the regularized field corresponds to a well defined measure $\mu_{k,V}$ by:
\begin{equation}
\label{Intro2}
R_{k,V}:=\frac{d\mu_{k,V}}{d\mu_{0}}= \frac{e^{-{\cal A}_{k,V}}}{Z_{k,V}} e^{-{\cal A}_{k,V}}, \; \;  
                              {\cal A}_{k,V}=\int_{V}L_{k}(\varphi_{k,V}(x), \partial \varphi_{k,V}(x))dx.
\end{equation}
\noindent
And the problem is whether a limit of $\mu_{k,V}$ exists in some sense, including that of the convergence of the Schwinger functions $S_{k,V}(f_{1}, ..., f_{p}) = \int \varphi(f_{1}) ... \varphi(f_{p}) d\mu_{k,V} $ and is non trivial. For the $\varphi^{4}_{d}$ field, let us write $L_{k}$ as (we drop the subscript $V$):
\begin{equation}
L_{k}(\varphi_{k,V}(x), \partial \varphi_{k,V}(x))= g_{k} :\!\varphi_{k}(x)^{4}\!: + m_{k} :\!\varphi_{k}(x)^{2}\!: + a_{k} :\!\partial \varphi_{k}(x)^{2}\!: - v_{n},
\end{equation}
where $g_{k}, m_{k}, a_{k}, v_{k}$ are renormalization constants of the coupling constant, mass, field strength and vacuum respectively (we use the notation $m_{k}$ instead of the usual $(\delta m_{k})^{2}$). For dimensions $d\geq 5$, a negative answer was obtained, in the case of one or two components $\varphi^{4}$ fields, by Aizenman \cite{Aizenman} and Fr\" ohlich \cite{Frohlich82}; using a lattice cutoff it was proved that the 
corresponding limiting field is Gaussian and hence a trivial one, for it is similar to the free field without interaction. We refer to Fernandez, Fr\" ohlich and Sokal \cite{FFS} for a detailed account on these questions and to Callaway \cite{Callaway} for a survey of the problem of triviality in QFT.
\\
The border case of dimension 4, which is the physical one, has so far remained an open problem, although partial results have been obtained for one and two component fields and with some conditions on the renormalization constants\footnote{For the $\varphi_{4}^{4}$ 1-component field, the remaining case corresponding to non vanishing $a_{k}$ has been treated recently by Aizenman and Dumenil-Copin \cite{A-DC} within the lattice approximation framework. See also the part II of this work \cite{AA3}.}; the case of negative coupling constant is also studied in Gawedzki, Kupiainen \cite{GK}.  It is believed that the limits of the regularized fields $\varphi_{k,V}(x)$ may also be trivial or that the interacting field may not exist, which would rise questions about the foundations and consistency of quantum field theory.
\\
\\
\noindent
The purpose of this paper is to show that for dimensions $d\geq 4$, and depending on the renormalization constants adopted, we have the following alternative: when the ultraviolet cutoff is removed and the volume is fixed: (1) either $\mu_{k,V}$ converges strongly to the Gaussian measure and $\varphi_{d}^{4}$ is trivial in this case, or (2) the almost sure limit of $R_{k,V}:= d\mu_{k,V} / d\mu_{0}$ exists and is equal to $0$. This implies, in the second case, which is the physical one, that $\mu_{k,V}$ can not have a strong limit as $k \longrightarrow \infty$. The second possibility is valid provided the coupling constant sequence of the modified action is not vanishing in the sense that $g_{k}k^{d-4}$ does not converge to $0$.
On the other hand, the proof of these results shows that the normalization sequence $Z_{k,V}=\E e^{-{\cal A}_{k,V}}$ is divergent as $k\rightarrow \infty$ when the vacuum renormalization $v_{k}$, is lower than some threshold, namely : $v_{k}\ll g_{k}c_{k}^{2}\sim g_{k}k^{2(d-2)} $, with $c_{k}:=\E \varphi_{k}(0)^{2}\sim k^{d-2}$. 
We recall that the boundedness of $Z_{k,V}$ is linked to the boundedness from below of the full Hamiltonian $H=H_{0}+{\it L}$ ({\it L} being the interaction Lagrangian). This was precisely the basic result of Nelson \cite{Nelson1} subsequently used to prove the existence of scalar field in dimension $2$, and in dimension $3$ this was the
main result of Glimm-Jaffe paper on the positivity of the $\varphi_{3}^{4}$ Hamiltonian \cite{GJ73}. The non positivity (or non boundedness from below) of the full interacting Hamiltonian in dimension 4 is an obstruction to the existence of non trivial interacting scalar quantum field in dimension 4 or greater.
\\
\\
The approach adopted in this paper is different from the previous ones. It uses direct calculations based on a momentum cutoff. A normalization of the field $\varphi$ and transformation of the integral of the Lagrangian to a mean of an array of random variables establish a link with the classical probabilistic questions of the law of large numbers and large deviations theory, and we are led to find an estimate of the repartition function $Z_{k,V}$ corresponding to $\mu_{k,V}$ via a Laplace type method.
\\
\\
\noindent
The remainder of the paper is organized as follows. The notations and estimates of the covariance functions are recalled in section 2. The statement of the result with some 
remarks are presented in section 3. Section 4 contains the proof, which is divided into four parts. The first part deals with a transformation of the Lagrangian and the action
to normalized ones and their expression as means of arrays of random variables; a law of large numbers and some estimates for this array of rv's are stated. 
In the second part, a lower bound of large deviations of the mentioned array and normalized Lagrangian is derived by using an approach of Bahadur, Zabell and Gupta \cite{BZG} ; we also use a result of Csisz\'ar \cite{Csiszar} that enables to identify the large deviations rate function. The third part uses this result to get a lower bound of the partition function $Z_{k,V}$ (by Varadhan's lemma). In the last part, the proof is completed with the different cases of the renormalization sequences. Finally, let us point out that the specific characters of the $\varphi^{4}$ or ferromagnetic fields does not play a role in the approach adopted in this work; the results obtained are also valid for multi-component (vector) fields and can be extended to polynomial Lagrangians.


\section{Notations and settings}

\subsection{Notations}

We consider the usual framework of Euclidean field theory. The probability space is $ {\cal P}= S'(\R^{d}), d \geq 1$
the Schwartz space of tempered distributions, with its Borel $\sigma$-field ${\cal B}$ ; the reference probability measure denoted by $\mu_{0}$ is the Gaussian measure whose covariance operator is $( -\Delta+1)^{-1}$. Unless otherwise specified, the norm $\|.\|_{p}$, will denote the $\L^{p}({\cal P})$ norm with respect to the measure $\mu_{0}$.
\\
We shall use a volume cutoff $V \subset \R^{d}$ and a momentum ultraviolet cutoff denoted by $n$ (instead of $k$):
\\
\noindent
For $k\geq 1$, let $\delta^{n} (x) = n^{d} \chi (n x )$ be a $ C^{\infty}$ approximation of the delta distribution ; then the 
regularized field is given by 
\begin{equation}
\label{RegField}
\varphi_{n}(x)=  \varphi ( \delta^{n}(. - x)),
\end{equation}
where $\varphi $ is the free field associated to $\mu_{0}$. The expectation with respect to  $\mu_{0}$ will be denoted by  $ \langle \rangle $ or $ \E $.
\\
\\
\noindent
If $P, Q$ are two probability measures on some space $E$, then $D(Q||P)$ will denote the Kullback-Leibler information (or number) or $I$-divergence of $Q$ with respect to $P$, that is:
\begin{equation}
\nonumber
 D(Q||P)=\left\{ \begin{array}{c}
 \int \log(dQ/dP)dQ \; \mbox{if} \; Q \ll P \\
 +\infty \; \; \mbox{otherwise}.
\end{array}\right.
\end{equation}


\subsection{Variances, covariances and estimates}

We shall use the following notations:
\begin{eqnarray*}
        c_{n}&= &  \langle \varphi_{n}(x)^{2} \rangle =   \langle \varphi_{n}(0)^{2} \rangle  \\
  c_{nl}(x,y)&= &  \langle \varphi_{n}(x) \varphi_{l}(y) \rangle =  \langle \varphi_{n}(x-y) \varphi_{l}(0) \rangle = c_{nl}(x-y), \\
      c_{n}(x-y)&= &  \langle \varphi_{n}(x) \varphi_{n}(y) \rangle = \langle \varphi_{n}(x-y) \varphi_{n}(0) \rangle  =c_{nn}(x-y)).
\end{eqnarray*} 
\noindent
NB. In the estimates we are concerned with, the constant factors will be denoted by the same letters $K, K'$, although they may be different and depend on the quantities estimated. In the case where these constants depend of some parameter e.g., $p$, this will be taken into account in the notation.
\begin{Rk}
The covariance $C(x,y)$ of the free field has the expression (see \cite{GJ87} pp. 162-163):
\begin{equation}
\label{covar1}
C(x,y) = \frac{1}{2(2\pi)^{d-1}} g(|x-y|), \; with \; g(t) = \frac{1}{t^{d-2}}|S^{d-2}|\int_{0}^{\infty} \frac{s^{d-3}e^{-s(1+t^{2}s^{-2})^{1/2}}}{(1+t^{2}s^{-2})^{1/2}} ds
\end{equation}
\noindent
From this expression it can be easily seen (see \cite{GJ87}) that:
\begin{equation}
\label{covar3}
g(t) \sim K t^{-(d-2)}
\end{equation}
Furthermore, if we take the derivatives of $g$ when $t > 0$, an inspection of the dominating terms in their expressions shows that we have also:
\begin{equation}
\label{covar4}
g'(t) \sim K -(d-2) t^{-(d-2)}, \; \; g''(t) \sim K (d-1)(d-2) t^{-(d-2)},
\end{equation}
and that in the manipulations involving the estimates of these terms, we are authorized to take only the previous dominating terms. We have for instance:
\begin{equation}
\label{covar5}
c_{n}(x) \sim  K \frac{1}{\sqrt{2}(2\pi)^{d/2}} \int_{\R^{d}} \frac{e^{-u^{2}/4}}{|u/n + x|^{d-2}} du
\end{equation}
\end{Rk}
As for the terms involved in the action, we set:
\begin{align*}
 I_{n}&=   \int_{V}  :\varphi_{n}(x)^{4}: dx ,\\
 M_{n}&=  \int_{V}  :\varphi_{n}(x)^{2}: dx ,\\
D_{n}&=   \int_{V}  : (\partial_{x}\varphi_{n}(x))^{2}:  dx ,
\end{align*} 
where $ :  : $ denotes the Wick product and $(\partial_{x}\varphi_{n}(x))^{2}$ is the norm of the vector 
$$(\partial_{x_{1}}\varphi_{n}(x), ..., \partial_{x_{d}}\varphi_{n}(x)).$$ 
We recall that if $f$ is a regular homogeneous random field, with covariance function $ \langle f(x) f(y) \rangle =c(x-y)$, then the variance of its derivatives is given by:
\[ \langle (\partial_{x_{j}} f(x))^{2} \rangle= -\frac{\partial^{2}c}{\partial x_{j}^{2}}(x) \]
We note that:
\begin{align*}
 \langle I_{n}^{2} \rangle&=   4! |V|\int_{V} c_{n}^{4}(x) dx,\\
 \langle M_{n}^{2} \rangle &=   2! |V| \int_{V} c_{n}^{2}(x) dx, \\
 \langle D_{n}^{2} \rangle &=   2! |V| \int_{V} \sum_{j=1}^{d}  ( \frac{\partial^{2} c_{n}}{\partial x_{j}^{2}}(x))^{2} dx
\end{align*} 
The covariances and their integrals depend on the dimension $d$, and we have the following estimates:
\begin{equation}
c_{n} \sim 
     \begin{cases}
          K \log n & \text{ if $d=2$}\\
         K n^{d-2} & \text{ if $d \geq 3$}
		\end{cases}	
\end{equation}

\begin{equation}
\int_{V} c_{n}^{2}(x) dx \sim
                        \begin{cases}
								           O(1) & \text{ if   $d = 3$} \\
                           K \log n &  \text{  if $ d=4$ }\\
                            K n^{d-4} & \text{ if $d \geq 5$}
											\end{cases}											
\end{equation}
\begin{equation}
\int_{V} c_{n}^{3}(x) dx   \sim   
															\begin{cases}
															K \log n &  \text{ if $  d=3$} \\
															K n^{2d-6} & \text{ if $ d \geq 4$ }
																\end{cases}																
\end{equation}
\begin{equation}
\int_{V} c_{n}^{4}(x) dx   \sim 
															\begin{cases}
																O(1) &  \text{ if  $d =2 $}  \\
																K n^{3d-8} & \text{ if   $d \geq 3$}
															\end{cases}
\end{equation} 
\noindent
As for the gradient field $\partial_{x}\varphi_{n}(x)$ we have the following result:
\begin{lm}
The gradient field $\partial_{x}\varphi_{n}(x)$ is a Gaussian field and its variance has the following estimate:
\begin{equation}
\label{vargrad}
\langle (\partial_{x}\varphi_{n}(x))^{2} \rangle= \langle (\partial_{x}\varphi_{n}(0))^{2} \rangle \sim K n^{d}
\end{equation}
\end{lm}
{\bf Proof.}\\
By a change of variable $v = u+nx$ in the estimate (\ref{covar5}) we have:
\begin{equation}
\label{covar6}
c_{n}(x) \sim  K n^{d-2} \frac{1}{\sqrt{2}(2\pi)^{d/2}} \int_{\R^{d}} \frac{e^{-(v-nx)^{2}/4}}{|v|^{d-2}} dv
\end{equation}
In view of Remark 2.1, by taking the derivatives:
\begin{equation}
\frac{\partial c_{n}}{\partial x_{j}}(x) \sim K n^{d-2} \frac{1}{\sqrt{2}(2\pi)^{d/2}} \int_{\R^{d}} \frac{2n (nx_{j}-v_{j}) e^{-(v-nx)^{2}/4} }{|v|^{d-2}} dv
\end{equation}
\noindent
and:
\begin{equation}
\label{SecondDerivative}
\frac{\partial^{2} c_{n}}{\partial x_{j}^{2}}(x) \sim K n^{d-2} \frac{1}{\sqrt{2}(2\pi)^{d/2}} [\int_{\R^{d}} \frac{4n^{2} (nx_{j}-v_{j})^{2} e^{-(v-nx)^{2}/4} }{|v|^{d-2}} dv
                                        + \int_{\R^{d}} \frac{ 2n^{2} e^{-(v-nx)^{2}/4} }{|v|^{d-2}} dv ]                                      
\end{equation}
\noindent
Which gives for $x=0$:
\begin{equation}
\frac{\partial^{2} c_{n}}{\partial x_{j}^{2}}(0) \sim 2 K n^{d} \frac{1}{\sqrt{2}(2\pi)^{d/2}} \int_{\R^{d}} \frac{(v_{j}^{2}+1) e^{-v^{2}/4} }{|v|^{d-2}} dv
                                                  \sim 2 K' n^{d}.
\end{equation}
This proves the lemma. $\Box$
\\
\noindent																				
Note that we can rewrite (\ref{SecondDerivative}) by remaking the change of variable $u = v-nx$:
\begin{equation}
\label{SecondDerivative2}
\frac{\partial^{2} c_{n}}{\partial x_{j}^{2}}(x) \sim 2K n^{d} \frac{1}{\sqrt{2}(2\pi)^{d/2}} \int_{\R^{d}} \frac{2(u_{j}^{2}+1) e^{-(u)^{2}/4} }{|u+nx|^{d-2}} du
\end{equation}
We shall use this estimate later.


\section{Statement of the results}

\subsection{The modified Lagrangian and random field}

With the notations of \S 2, we consider the Euclidean quantum field with interaction, which is associated to a Lagrangian density $L(\varphi(x))$ and the action:

\begin{equation}
\label{lagrangian}
{\cal A} = \int_{R^{d}}L(\varphi(x))dx  
\end{equation}
\noindent
and to a measure $\mu$ on ${\cal P}$ given by:
\begin{equation}
\label{field1}
\frac{d\mu}{d\mu_{0}}= \frac{1}{Z}\exp[-\int_{R^{d}}L(\varphi(x)dx)], \; \;  Z =E_{\mu_{0}}\exp[-\int_{R^{d}}L(\varphi(x)dx)] 
\end{equation}
In the scalar field case (bosonic interactions), the Lagrangian density $ L$ is usually a polynomial function $P(\varphi(x))$, and the most simple but fundamental scalar interaction is the case where $P(\varphi(x))= \varphi(x)^{4}$ or $P(\varphi(x))=:\varphi(x)^{4}:$. Yet, the expressions (\ref{lagrangian}) and (\ref{field1}) are formal because $\varphi$ is a distribution. So regularizations are performed via several methods and we shall use the momentum and volume regularizations recalled 
in $\S 2$ : $\varphi$ will be replaced by a regular field $\varphi_{n}$ and $\R^{d}$ by a finite volume $V$.
\\
On the other hand, in the study of the limit of the interacting regularize field, we are led to add counter-terms to the the interactions in order to get 
meaningful and finite quantities: in the perturbative (or physical) renormalization, the quantities in question are the moments of the type $\langle \varphi(x_{1}) \varphi(x_{2}) ... \varphi(x_{p})\rangle$. In the constructive renormalization, the aim is to obtain the limit of the field itself as a well defined and non trivial object that gives finite moments, and to see whether they are the same as the moments obtained by the perturbative procedure (cf., e.g., Rivasseau \cite{Rivasseau}). The outcome depends on the space-time dimension $d$.
\\
\\
\noindent
$\bullet$ {\it Dimension 2:}  The Wick regularization, replacing $\varphi(x)^{4}$ by $:\varphi(x)^{4}:$, is sufficient to construct an interacting field in infinite volume ; this case had been studied by Nelson \cite{Nelson1}, Glimm, Jaffe \cite{GJ68} and many authors by the end of the 1960s and beginning of the 1970s, cf. Glimm-Jaffe \cite{GJ87} and Simon \cite{Simon} for references. 
\\
$\bullet$ {\it Dimension 3:} The Wick regularization is not sufficient: a mass renormalization term is added and the modified Lagrangian has the form:
\begin{equation}
\label{RlagrangianD3}
L_{n} (\varphi) = g  :\varphi_{n}^{4}(x): +  (\delta m_{n})^{2}  :\varphi_{n}^{2}(x): dx 
\end{equation}
\noindent
In both cases (d=2,3), the renormalization scheme is similar to that of the perturbative renormalization, cf. Gallavotti \cite{Gallavotti1}, \cite{Gallavotti2}. Still, the proof of the existence of these fields {\it and} their non triviality is a highly non-trivial task and is considered as a major achievement in the constructive quantum field theory program. In fact the difficulty is rather
to prove that a regularized field with a given renormalization constants converges to some field that is non trivial and satisfies the OS axioms.
\\
\\
\noindent
$\bullet$ {\it Dimension 4:} In this case, besides the mass renormalization term $(\delta m_{n})^{2}$, the perturbative renormalization procedure requires the addition of other counter-terms, namely the constant coupling and wave function renormalization (see e.g \cite{Kleinert-S}, \cite{Collins}, \cite{Salmhofer}); the modified Lagrangian density will have the following form:
\begin{align}
\label{RlagrangianD4}
L_{n} (\varphi) &=  g_{n} :\varphi_{n}^{4}(x): + (\delta m_{n})^{2}  :\varphi_{n}^{2}(x):  + a_{n} :(\partial \varphi_{n})^{2}(x): \nonumber \\
                       &= L_{n}(\varphi_{n}(x), \partial \varphi_{n}(x)) dx
\end{align}
{\it Remark and notations:} In the complete Lagrangian density we have to add a vacuum renormalization term $v_{n}$, so that the full Lagrangian density is $L_{n}^{*}=L_{n}-v_{n}$:
\begin{equation}
\label{RlagrangianD4b}
L_{n}^{*}(\varphi_{n}(x), \partial \varphi_{n}(x)) =  g_{n} :\varphi_{n}^{4}(x): + (\delta m_{n})^{2} :\varphi_{n}^{2}(x):  + a_{n} :(\partial \varphi_{n})^{2}(x): -v_{n}         
\end{equation}
and the complete action when we consider a finite volume $V$ will be:
\begin{equation}
\label{RlagrangianD4b}
{\cal A}_{n}^{*}=\int_{V}g_{n} [:\varphi_{n}^{4}(x): + (\delta m_{n})^{2} :\varphi_{n}^{2}(x):  + a_{n} :(\partial \varphi_{n})^{2}(x): -v_{n}] dx 
\end{equation}
Except some remarks below, this vacuum renormalization term will be left aside in the sequel: the interacting field measures given by $d\mu_{n}=e^{-{\cal A}_{n}^{*}}/Z^{*}_{n}d\mu_{0}$, with $Z^{*}_{n}= \E e^{-{\cal A}_{n}^{*}} $ are the same as those given by $d\mu_{n}=e^{-{\cal A}_{n}}/Z^{*}_{n}d\mu_{0}$, with $Z_{n}= \E e^{-{\cal A}_{n}} $, the term $e^{-v_{n}}$ being eliminated between $e^{-{\cal A}_{n}^{*}}$ and $Z^{*}_{n}$ in the first expression.

In this paper we shall consider the following action :
\begin{equation}
\label{Action}
{\cal A}_{n}=\int_{V}g_{n} [:\varphi_{n}^{4}(x): +  m_{n} :\varphi_{n}^{2}(x):  + a_{n} :(\partial \varphi_{n})^{2}(x):]dx
\end{equation}
where $V\subset \R^{d}$ is a finite volume. The mass term renormalization will be denoted by $m_{n}$ instead of the usual $(\delta m_{n})^{2}$. We call $g_{n}, m_{n}, a_{n}$ the renormalization constant sequences or simply the renormalization constants.
In order to simplify the notations we will often set: ${\cal A}_{n}(\varphi) =  \int_{V} L_{n}(\varphi_{n}(x)) dx$ and it is understood that the Lagrangian depends also on the gradient of the field. The interacting field measure $\mu_{n,V}$ for $\varphi_{d}^{4}$ in a finite volume and with a momentum cutoff $n$ is given by:
\begin{equation}
\label{IFmeasure}
\frac{d\mu_{n,V}}{d\mu_{0}}= \frac{1}{Z_{n,V}}e^{ -{\cal A}_{n}}, \; \; Z_{n,V}= \E e^{ -{\cal A}_{n}}
\end{equation}

\subsection{Main result and remarks}

The following theorem and corollaries are the main result of this paper:

\begin{thm}\label{main-th}
Let $\mu_{n,V}$ be the probability measure on ${\cal S}'(\R^{d})$ associated to the $\varphi_{d}^{4}$ field with a momentum and a volume cutoffs indexed by $n, V$, and
renormalization constant sequences $g_{n}, m_{n}$ and $a_{n}$. We suppose that the coupling constant sequence is positive $g_{n} \geq 0$. Then, with $c_{n}= \E\varphi_{n}(0)^{2}\sim K_{d} n^{d-2}$ and $d_{n}\sim n^{2}$, we have the following alternative:
\\
(1) If the sequences $g_{n} c_{n}^{2}, m_{n}c_{n}$ and $a_{n} d_{n }c_{n}$ are bounded, then $\mu_{n,V}$ converges strongly (setwise) to the Gaussian measure $\mu_{0}$ and $\varphi_{d}^{4}$ is therefore the free field.
\\
or
\\
(2) If at least one of the sequences $g_{n} c_{n}^{2}, m_{n}c_{n}$ or $a_{n} d_{n }c_{n}$ is unbounded then, provided that the dimension $d\geq 4$ and either the coupling constant sequence is not vanishing in the sense that $g_{n}n^{d-4}$ does not converge to $0$, or if this condition is not fulfilled: $m_{n}\geq K' n^{2}$ or $a_{n}\geq K'$ for some constant $K' >0$, we have:
\begin{equation}
\label{mainth1}
\lim_{n \rightarrow \infty} \frac{d\mu_{n,V}}{d\mu_{0}}= \lim_{n \rightarrow \infty}\frac{1}{Z_{n,V}}\exp [-{\cal A}_{n} ]= 0 \; \;a.e,  
\end{equation}
where 
\begin{equation}
\label{mainth2}
Z_{n,V}=E_{\mu_{0}} \exp [-{\cal A}_{n} ]
\end{equation}
and in this case there exists a constant $K >0$ such that:
\begin{equation}
\label{mainth3}
Z_{n,V} \geq e^{K \max (g_{n} c_{n}^{2}, m_{n}c_{n}, a_{n} d_{n }c_{n})}
\end{equation}
The previous assertions are also valid in the case the $|\varphi|^{4}$ multi-component field and for all dimensions $d\geq 4$.
\end{thm}

\noindent
The following first corollary shows that in dimensions $d\geq 4$, it is not possible to construct an interacting $\varphi^{4}$ field in a finite volume $V$ as a strong limit of its cutoff measure $\mu_{n,V}$:

\begin{co}\label{main-coro}
With the notations of the previous theorem, in the case (2) where one of the sequences $g_{n} c_{n}^{2}, m_{n}c_{n}$ or $a_{n} d_{n }c_{n}$ is unbounded with one of the further conditions mentioned in the theorem, the sequence:
$$ R_{n,V}= \frac{d\mu_{n,V}}{d\mu_{0}} $$
is not uniformly integrable and the sequence of measures $\mu_{n,V}$ does not converge strongly (setwise) to any probability measure. 
\end{co}

\noindent
{\bf Proof.}\\
If  $R_{n,V}$ is uniformly integrable, then, since $R_{n,V}  \longrightarrow R_{\infty}=0 $ a.e,  we will have $ \E R_{n,V}  \longrightarrow \E R_{\infty}=0 $ (by the extension 
of the dominated convergence theorem), which is not possible because $\E R_{n,V}=1$.
\\
\noindent
Now, if  $\mu_{n,V}$ converges setwise, i.e. $\mu_{n,V} (A)$ converges to some $\mu_{\infty} (A)$ for every $A\in {\cal B}$, then $\int_{A} R_{n,V} d \mu_{0}$  for every $A\in {\cal B}$. But this would imply that $R_{n,V}$ is uniformly integrable by the Vitali-Hahn-Sacks theorem (cf., e.g., Neveu \cite{Neveu}, Proposition IV-2.2), which is in
contradiction with the first part of the corollary. $\Box$.
\\
\\
\noindent
The second corollary concerns the ultraviolet stability of the $\varphi_{d}^{4}$ model:
\begin{co}\label{main-coro2}
Let ${\cal A}_{n}^{*}$ be the complete renormalized action of the $\varphi_{d}^{4}$ model in a finite volume $V$ and dimension $d\geq 4$:
\begin{equation}
\label{coro2a}
{\cal A}_{n}^{*}=\int_{V}g_{n} :\varphi_{n}^{4}(x): + \delta m_{n} :\varphi_{n}^{2}(x):  + a_{n} :(\partial \varphi_{n})^{2}(x): -v_{n} 
\end{equation}
 and $Z^{*}_{n}= \E e^{-{\cal A}_{n}^{*}}$. We also assume that the coupling constant $g_{n}$is not vanishing. Then for any renormalization scheme in which the vacuum renormalization $v_{n}$ satisfies $v_{n} \ll g_{n} n^{2d-4}$, we will have:
\begin{equation}
\label{coro2b}
Z^{*}_{n}= \E e^{-{\cal A}_{n}^{*}} \rightarrow +\infty
\end{equation}
The corresponding Hamiltonian of the limiting field will not be bounded from below (and can not provide a quantum field theory model)
\end{co}
\noindent
{\bf Proof.}\\
Since $g_{n} c_{n}^{2} \sim K g_{n} n^{2d-4} $ and $\limsup g_{n} > 0$, we will have in the case where $v_{n} \ll g_{n} n^{2d-4}$ : 
$$\E e^{-{\cal A}_{n}^{*}}  \geq e^{Kg_{n} c_{n}^{2}-v_{n}}\geq  e^{K'g_{n} n^{2d-4} -v_{n}} \rightarrow +\infty $$
for some constant $K' >0$. For the second assertion we refer to the remarks below. $\Box$
\\
Let us now give some comments on the these results:
\\
\\
1. {\it On the role of the space-time dimension:} although the possibility $(1)$ of the theorem is valid for all dimensions $d\geq1$, there is no interference between theorem 3.1 and the non-trivial constructions of quantum field carried out for dimensions $d=2,3$:
\\
\\
\noindent
$\bullet$ In dimension $d=2$, we have $c_{n} \sim K \log n$ and the construction of of quantum fields in this case uses $g_{n}= Constant$ or
$O(1)$, while there is no need of mass and wave renormalizations. Theorem 1 does not apply to this configuration.
\\
\\
\noindent
$\bullet$ In dimension $d=3$, we have $c_{n} \sim K n$, and the construction of of quantum fields in this case uses $g_{n}= O(1)$ and $m_{n} \sim K \log n$. Hence, the term $g_{n} c_{n}^{2}$ is unbounded, so we are in the Case 2 of the theorem which requires $d\geq 4$. In fact, we will see in the proof that the requirement $d\geq 4$ is set in order to fulfill the condition : 
\begin{equation}
\label{condition}
g_{n} c_{n}^{2} \sim g_{n} n^{2d-4}\geq K n^{d},
\end{equation} 
which is satisfied in dimensions $d\geq 5 $ and $d= 4 $ if $g_{n}\nrightarrow 0$. (\ref{condition}) would be satisfied in dimension 3 if the coupling constant is divergent and $g_{n}=O(n^{2})$, which is not the renormalization scheme adopted in this dimension.
\\
\\
\noindent
3. In the constructive quantum field theory literature the uniform boundedness of $Z_{n,V}$ (sometimes together with the existence of a limit field as the UV cutoff is removed) 
is usually referred to as the ultraviolet stability of the model. That $Z_{n,V}:= \E \exp(-{\cal A}_{n}^{*})$ is uniformly bounded (in $n$, the volume $V$ being fixed) is indeed a key step towards the construction of mathematical QFT models. In dimension 2, this result was obtained by Nelson \cite{Nelson1} and proved in other ways by many authors, 
see \cite{Simon} for references. It was subsequently used by Glimm and Jaffe for the first construction of scalar quantum field in dimension $2$ in finite volume within the Hamiltonian framework. For dimension 3, the uniform bound for $Z_{n,V}$ was the main purpose of the Glimm-Jaffe article, Positivity of the $\varphi_{3}^{4}$ Hamiltonian \cite{GJ73}, with its consequence, the uniform boundedness from below of the total interacting Hamiltonian $H_{n,V}= H_{0}+{\cal L}_{n,V}$. Besides being a physical requirement, the boundedness from below of the total Hamiltonian is used to prove the essential self-adjointness of the Hamiltonian $H_{V}$ obtained as the UV cutoff is removed, which enables to continue the construction of the interacting quantum field. The other consequence of this positivity of the Hamiltonian is that $\inf {\rm spectrum} (H) > -\infty$ which is a condition for the existence of the lowest energy state, i.e. the vacuum.
The Feynman-Kac-Nelson formula (used for fields with ultraviolet cutoff) provides a link between the boundedness from below of the total Hamiltonian and the ultraviolet stability (see \cite{GJ73}, \cite{Simon}, \cite{GRS75}). The divergence of $Z_{n,V}^{*}$ stated above (besides the trivial cases) presents therefore a problem for the construction of a scalar quantum field and its possible existence in dimensions $d\geq 4$ (if $v_{n} \ll g_{n} n^{2d-4}$).
\\
\\
\noindent
4. We recall that a standard method to prove the existence of a field in $d$-dimension, is to show that the Schwinger functions:

\[ S_{n,V}(f_{1}, ..., f_{p}) = \int \varphi(f_{1}) ... \varphi(f_{p}) d\mu_{n,V} \]

\noindent
or  $ S_{\delta,V} (x_{1}, ..., x_{p})$ in the case of lattice regularization, have a limit as the cutoff are removed and to prove the properties related to the Osterwalder-Schrader axioms. One can also seek the limits of the characteristic functionals  $C_{n,V}(f) = \int_{V} \exp (i\varphi(f)) d \mu_{n,V}$. We also recall that when taking the infinite volume limit we can not hope a strong convergence of the measures $ \mu_{n,V}$ to a non trivial measure $ \mu_{\infty} $: we would then have $ \mu_{\infty} \ll \mu_{0} $ but due to the Haag theorem (Euclidean version, see \cite{Fr74}, \cite{Schrader1}, \cite{LN}), $\mu_{\infty}$ would be a Gaussian measure.
\\
However, in finite volume it is possible to have a non trivial strong limit of $ \mu_{n,V}$ (The Haag theorem is valid only in infinite volume). This has been accomplished in dimension 2, see, eg. Newman \cite{Newman1}, and the remarks in Simon \cite{Simon}, p. 142.
\\
\\
\noindent
In dimensions $d\geq 4$, this possibility of a non trivial strong limit of $ \mu_{n,V}$ with a fixed finite volume is nevertheless ruled out by theorem \ref{main-th} and corollary \ref{main-coro}
\\
Several attempts have been made in 4 dimension, see e.g. \cite{Schrader}, \cite{GR84}. In relation with, this point Glimm and Jaffe made a remark in \cite{GJ74-2}, that for a class of $\varphi^{4}$ fields, a bound of the two point function of the type $|S_{\varepsilon}^{(2)}(f\otimes g)| \leq |f|_{\cal S}|g|_{\cal S}$ is sufficient to prove a bound on the n-point Schwinger function $|S_{\varepsilon}^{(n)}(x_{1}, ..., x_{n})|$ independent of $\varepsilon$, in the lattice framework, which yields the existence of  $|S^{(n)}(x_{1}, ..., x_{n})|$ as a limit when the lattice cutoffs $\varepsilon \longrightarrow 0$, by a compactness argument. One has to prove then the OS axioms. Still, the question of non-triviality of the field obtained has to be addressed and may be more difficult.
\\
\\
\noindent
5. For a discussion about the triviality concept we refer to \cite{GR84}, \cite{FFS} and \cite{Callaway}. A trivial limit includes the cases of free field (no interaction) or singular field; the limit of the regularized field (or measure) may also not exist at all. The standard method to prove that the limiting field is trivial within the lattice framework is to show that the Ursell functions (truncated 4-point function):
\begin{align}
\label{Ursell}
U_{4}(\varepsilon):= \langle\varphi^{i_{1}} \varphi^{i_{2}} \varphi^{i_{3}} \varphi^{i_{4}}\rangle^{T}_{\varepsilon} & :=  
                                      \langle\varphi^{i_{1}} \varphi^{i_{2}} \varphi^{i_{3}} \varphi^{i_{4}}\rangle_{\varepsilon} \\ \nonumber
                           & - \langle\varphi^{i_{1}} \varphi^{i_{2}}\rangle_{\varepsilon} \langle\varphi^{i_{3}} \varphi^{i_{4}}\rangle_{\varepsilon} 
													- \langle\varphi^{i_{1}} \varphi^{i_{4}}\rangle\langle\varphi^{i_{2}} \varphi^{i_{3}}\rangle_{\varepsilon}
\end{align}
converges to $0$ as the the lattice spacing $\varepsilon$ goes to $0$. In this expression $\langle \rangle_{\varepsilon}$ denotes the expectation w.r.t a lattice measure $\mu_{\varepsilon}$ approximating the scalar field. That $U_{4}(\varepsilon)\rightarrow 0$ implies that the resulting field is trivial, is a consequence of either a theorem of Baumann \cite{Baumann} or a theorem of Newman \cite{Newman75}, which assert that when $U_{4}=0$, the field under consideration is a generalized free field, see \cite{FFS} for the validity conditions of these theorems and other details.
\\
\\
\noindent
6. {\it The case $g_{n}\rightarrow 0$:} In the lattice framework, the limiting field is trivial in this case and this is a consequence of the skeleton inequalities (\cite{BFS83}) :
\begin{equation}
\label{skeleton}
0 \leq -U_{4}  \leq 4! \lambda \sum_{j\in \L} \langle\varphi^{i_{1}} \varphi^{j}\rangle\langle\varphi^{i_{2}} \varphi^{j}\rangle\langle\varphi^{i_{3}} \varphi^{j}\rangle\langle\varphi^{i_{4}} \varphi^{j}\rangle 
\end{equation}
(the correspondence between coupling constant in the lattice ($\lambda$) and continuum framework ($g$) is: $\lambda=\varepsilon^{d} g$). In our case we have the following situations :
\\
$\bullet$ Case (a) : $g_{n}\rightarrow 0$ and $g_{n}n^{d-4}\nrightarrow 0$: this happens in dimensions $d\geq 5$: in this case the results of theorem 3.1 are valid.
\\
$\bullet$ Case (b) : $g_{n}n^{d-4}\rightarrow 0$: theorem 3.1 does no apply in this case; that the resulting field is also trivial in this case does not seem to be straightforward from the arguments of the proof below, besides the case (1) of the previous theorem, where $g_{n}$ is rapidly vanishing ($g_{n}n^{2d-4}\rightarrow 0$). However it can be shown from the results of the lattice framework (the skeleton inequalities), that the limiting field is actually trivial. 

We leave aside the details for the moment. Let us just mention that in the two situations (a) and (b) the limiting field might be null, in particular when $m_{n}\gg g_{n}n^{d-2}$ (cf. \cite{AA2}).
\\
\\
\noindent
7. $P(\varphi)_{4}$: The proof below is valid without modification of the arguments to cover the case of polynomial interaction of course with the condition the polynomial $P$ is bounded below.


\section{Proofs}

\subsection{Overview of the proof}

The idea of the proof is to show that $Z_{n,V}$ is much larger that the values of:
$$ \exp [-\int_{V}L(\varphi_{n}(x))dx] \; {\rm as } \;  n\longrightarrow \infty $$ 
We therefore seek an estimate (in fact a minorization) of $Z_{n,V}$ and this is done via the following steps:
\\
$\bullet$ (1) $Z_{n,V}$ is first transformed to an expression of the form:
\begin{equation}
\label{overview1}
Z_{N,V} = \int_{B} e^{-N F(\varphi)} \mu_{N} (d\varphi), \; \; N=n^{d}
\end{equation}
\\
\noindent
$\bullet$ (2) Estimating expressions like (\ref{overview1}) is usually done with a Laplace type method. In infinite dimension, the Varadhan lemma, which is used in many situations, requires that the sequence of measures $\mu_{N}$ satisfies a large deviations principle (LDP), with assumptions that can be more or less relaxed depending of the situations. However, proving a LDP in our case seems to be complicated.
\\
$\bullet$ (3) Nonetheless, it turns out that a lower bound of large deviations for the sequence $\mu_{N}$ can be established and this is sufficient to get a minorization of $Z_{n,V}$. 
\\
$\bullet$ (4) To get this lower bound, the action and Lagrangian are transformed to arrays of random variables with an expression like: 
\begin{equation}
\label{overview2}
{\cal A}_{n} = u_{n} {\cal A}_{n}^{(1)}, \; \; {\cal A}_{n}^{(1)}=\int_{V}l_{n}(\varphi_{n}(x)dx = \frac{1}{n^{d}} \sum_{1}^{n^{d}} X_{n,i} 
\end{equation}
\\
\noindent
We are thus led to the study of the large deviations properties of the array $X_{n,i}$ and this is the major part of the proof.

\subsection{Transformation to an array of random variables and a law of large numbers for the integrated fields}

We normalize the field $\varphi_{n}$ by the field  $\psi$ defined by: 
$$       \psi_{n}(x)= \frac{\varphi_{n}(x)}{\sqrt{c_{n}}}  $$
\noindent
The renormalized action $ {\cal A}_{n}$ may be written as:
\begin{align*}
 {\cal A}_{n} &=  \int_{V}L(\varphi_{n}(x))dx  \\
              &= g_{n} \int_{V} :\varphi_{n}^{4}(x): dx + m_{n} \int_{V} :\varphi_{n}^{2}(x): dx + a_{n} \int_{V}:(\partial \varphi_{n})^{2}(x): dx  \\
              &=  g_{n} c_{n}^{2} \int_{V} :\psi_{n}^{4}(x): dx + m_{n}c_{n}\int_{V} :\psi_{n}^{2}(x): dx + a_{n}c_{n} d_{n} \frac {1}{d_{n}} \int_{V}:(\partial \psi_{n})^{2}(x): dx
\end{align*}
where we have added a factor $d_{n}$ in the last term that will be justified later. 
Depending on which of the factors $g_{n} c_{n}^{2}$, $m_{n}c_{n}$ or $a_{n}c_{n} d_{n}$ is the dominant one, we shall define a new Lagrangian $l_{n}$ density and action ${\cal A}_{n}^{(1)}$ by factorizing by the dominant term. 
For example when $g_{n} c_{n}^{2} \geq K m_{n}c_{n} $ $g_{n} c_{n}^{2} \geq K' a_{n}c_{n} d_{n} $, with $K, K' > 0 $ we write :
\begin{equation}
\label{ReduLag0}
L(\varphi_{n}(x))dx= g_{n} c_{n}^{2}l_{n}(\varphi_{n}(x))
\end{equation}
\[  l_{n}(\varphi_{n}(x)) = \lambda_{n} :\psi_{n}^{4}(x): + \alpha_{n} :\psi_{n}^{2}(x): + \beta_{n}\frac {1}{d_{n}} (\partial \psi_{n})^{2}(x):  \]
\begin{equation}
\label{ReduLag}
{\cal A}_{n}^{(1)}(\varphi)=   \lambda_{n}\int_{V} :\psi_{n}^{4}(x): dx + \alpha_{n}\int_{V} :\psi_{n}^{2}(x): dx + \beta_{n}\frac {1}{d_{n}} \int_{V}:(\partial \psi_{n})^{2}(x): dx
\end{equation}
and we have in this case: $\lambda_{n}=1$, $\alpha_{n}=m_{n}/(g_{n}c_{n})$ and $\beta_{n}=a_{n}d_{n}/(g_{n}c_{n})$.
\\
\\
\noindent
{\bf Notations related to the volume V. }
\\
For definiteness, the finite volume $V$ will be taken as $V=[0,1]^{d}$ and for each $n$, the volume $[0,n]^{d}$ will be denoted by $V^{(n)}$ and will be divided to a subdivision of $n^{d}$ volumes $V^{(n)}_{i}= [i_{1},i_{1}+1]\times [i_{2},i_{2}+1] .... \times [i_{d},i_{d}+1] $, with $ i= (i_{1}, i_{2}, ...., i_{d})$ and $i_{k}=0, ..., n-1$. There $  n^{d}$ indices $i$, and for the convenience of the notations used in the summations we shall make this index $i$ running in the set $\{1, ...,  n^{d}\}$ instead of the set $ i= (i_{1}, i_{2}, ...., i_{d}), i_{k}=0, ..., n-1 $ by a trivial correspondence. We also set:

\[ \frac{V^{(n)}_{i}}{n} = [\frac{i_{1}}{n},\frac{i_{1}+1}{n}]\times [\frac{i_{2}}{n},\frac{i_{2}+1}{n}] .... \times [\frac{i_{d}}{n},\frac{i_{d}+1}{n}]  \]

\noindent
The small volumes $V^{(n)}_{i}/n$ form a subdivision  of $[0,1]^{d} = V = V^{(n)}_{0}$. 

\noindent
Let us make the following transformation:
$$ X_{n,i} = \int_{V^{(n)}_{i}} l_{n}(\psi_{n}(\frac{x}{n})) dx  = n^{d} \int_{V^{(n)}_{i}/n} l_{n}(\psi_{n}(x)) dx $$

\noindent
${\cal A}_{n}$ and ${\cal A}_{n}^{(1)}$  are thus expressed as a mean of an array of the random variables $ X_{n,i}$:
\begin{equation}
\label{array1}
{\cal A}_{n}^{(1)}(\varphi)= \int_{V} l_{n}(\psi_{n}(x)) dx = \sum_{i=1}^{n^{d}}  \int_{V^{(n)}_{i}/n} l_{n}(\psi_{n}(x)) dx = \frac{1}{n^{d}}\sum_{i=1}^{n^{d}} X_{n,i}
\end{equation}
\noindent
or by a change of variable $y=nx$:
\begin{equation}
\label{array2}
{\cal A}_{n}^{(1)}(\varphi) = \frac{1}{n^{d}} \int_{V^{(n)}} l_{n}(\psi_{n}(\frac{y}{n})) dy = \frac{1}{n^{d}} \sum_{i=1}^{n^{d}}  \int_{V^{(n)}_{i}} l_{n}(\psi_{n}(\frac{y}{n})) dy= \frac{1}{n^{d}}\sum_{i=1}^{n^{d}} X_{n,i}
\end{equation}
\noindent
Let us introduce the following notations:
\begin{equation}
\label{IMD1}
I_{n,i}= \int_{V^{(n)}_{i}} :\psi_{n}^{4}(\frac{x}{n}): dx, \; \;  M_{n,i}= \int_{V^{(n)}_{i}} :\psi_{n}^{2}(\frac{x}{n}): dx,  \; \;  D_{n,i}=  \frac {1}{d_{n}}\int_{V^{(n)}_{i}}  :(\partial \psi_{n}(\frac{x}{n}))^{2}:dx
\end{equation}

\begin{equation}
\label{IMD2}
I_{n}= \int_{V} :\psi_{n}^{4}(x): dx, \; \;  M_{n}= \int_{V} :\psi_{n}^{2}(x): dx,  \; \;  D_{n}=  \frac {1}{d_{n}}\int_{V}  :(\partial \psi_{n}(x)^{2}:dx
\end{equation}

\noindent
The following proposition shows that the array of the random variables $ X_{n,i}$ has correct properties like the convergence of the $\L_{p}$ norms of the elements $ X_{n,i}$ to
non trivial values (not null and not infinite) and that they are asymptotically decorrelated.
\begin{prop}\label{PropArrays}
The arrays of random variables $I_{n,i}, M_{n,i}, D_{n,i}, X_{n,i}$ have the following properties:
\\
$\bullet$ (1) There exists constants $K_{p}, K_{D} $ such that, with $d_{n}=n^{2}$, we have the following limits:  
\begin{align*}
\label{prop14}
 \lim_{n\longrightarrow \infty} \langle (\int_{V^{(n)}_{i}} :\psi_{n}^{p}(\frac{x}{n}): dx)^{2} \rangle &= K_{p} \\
 \lim_{n\longrightarrow \infty} \langle (\frac{1}{d_{n}}\int_{V^{(n)}_{i}}  :(\partial \psi_{n}(\frac{x}{n}))^{2}:)^{2} \rangle &= K_{D}
\end{align*}

$\bullet$ (2) In particular: $\langle I_{n,i}^{2}\rangle \sim K_{4} $ and $ \langle M_{n,i}^{2} \rangle \sim K_{2}$ and if the sequences $\lambda_{n}, \alpha_{n}, \beta_{n}$ converge to $\lambda, \alpha, \beta$ then the $L_{2}$-norm of the modified array:
\begin{equation}
X_{n,i} = \lambda_{n} \int_{V^{(n)}_{i}} :\psi_{n}^{4}(\frac{x}{n}): dx + \alpha_{n}\int_{V^{(n)}_{i}} :\psi_{n}^{2}(\frac{x}{n}): dx + \beta_{n}  \frac {1}{d_{n}} \int_{V^{(n)}_{i}}  :(\partial \psi_{n})^{2}(\frac{x}{n}): dx
\end{equation}
converges to a real $K$ independent of $i$, that is $\lim_{n\longrightarrow \infty} \| X_{n,i}\|_{2}=K$.
\\
$\bullet$ (3) For $p\in\N$, the $\L_{p}$-norms $\| I_{n,i}\|_{p}, \| M_{n,i}\|_{p}, \| D_{n,i}\|_{p}, \| X_{n,i}\|_{p}$ are also convergent and hence bounded.
\\
$\bullet$ (4) For all $ i, j$ with $i \neq j$, $I_{n,i}$ and  $I_{n,j}$ are asymptotically decorrelated:
$$ \lim_{n\longrightarrow \infty} \langle I_{n,i} I_{n,j} \rangle =0, $$ 
\noindent
and the same property holds for $M_{n,i}, D_{n,i}, X_{n,i}$.
\end{prop}

\noindent
{\bf Proof.}
\\
To begin with (1), we have:
\begin{align*}
<I_{n,i}^{2}\rangle &= n^{2d}\langle (\int_{V^{(n)}_{i}/n} \int_{V^{(n)}_{i}/n} :\psi_{n}^{p}(x)::\psi_{n}^{p}(y): dxdy \rangle \\
              &=  n^{2d}\langle (\int_{V^{(0)}_{i}/n} \int_{V^{(n)}_{0}/n} :\psi_{n}^{p}(x)::\psi_{n}^{p}(0): dxdy \rangle \\
              &=  \frac{n^{2d}}{2^{d}}  \int_{ 2 V^{(n)}_{0}/n} \int_{2 V^{(n)}_{0}/n} \langle :\psi_{n}^{p}(x)::\psi_{n}^{p}(0):\rangle dxdy 
\end{align*}
where we have made a first change of variable $x\longrightarrow x -i/n$, $y\longrightarrow y -i/n$ with $i= (i_{1}, i_{2}, ...., i_{d})$ and a second one: 
$x\longrightarrow x +y $, $y\longrightarrow x-y $. We also set $2 V^{(n)}_{0}/n= [0, 2/n]^{d}$ with its volume $ |2 V^{(n)}_{0}/n|= (2/n)^{d}$. This yields:
\begin{align*}
     <I_{n,i}^{2}\rangle  &=  \frac{n^{2d}}{2^{d}}  \frac{2^{d}}{n^{d}} \int_{ 2 V^{(n)}_{0}/n} p! \langle \psi_{n}(x)\psi_{n}(0)\rangle^{p} dx \\
		                &=  \frac{n^{2d}}{2^{d}}  \frac{2^{d}}{n^{d}} \frac{p!}{c_{n}^{p}} \int_{ 2 V^{(n)}_{0}/n}  \langle \varphi_{n}(x)\varphi_{n}(0)\rangle^{p} dx \\
	                  &=  \frac {n^{d} p!}{(c_{n})^{p}} \frac{p!}{c_{n}^{p}} \int_{ 2 V^{(n)}_{0}/n}  c_{n}^{p}(x) dx \\	   
									  & \sim  \frac{p! n^{d}}{(c_{n})^{p}} \int_{2 V^{(n)}_{0}/n} dx \int_{\R^{d}} ... \int_{\R^{d}} 
											       \frac{1}{(\sqrt{2}2\pi)^{dp/2}}  \frac{e^{-(\sum_{j=1}^{p}u^{2}_{j})/2}}{\prod_{j=1}^{p}|u_{j}/n + x|^{d-2}} du_{1} ... du_{p} \\
										& \sim  \frac{p! n^{d} n^{-d} n^{p(d-2)}}{c_{n}^{p}}  \frac{1}{(\sqrt{2}2\pi)^{dp/2}}  \int_{2 V^{(n)}_{0}} dx' \int_{\R^{d}} ... \int_{\R^{d}}
										           \frac{e^{-(\sum_{j=1}^{p}u^{2}_{j})/4}}{\prod_{j=1}^{p}|u_{j} + x'|^{d-2}} du_{1} ... du_{p}, \\			
									  & \sim   \frac{p!}{(\sqrt{2}2\pi)^{dp/2}} \int_{\R^{d}} dx' \int_{\R^{d}} ... \int_{\R^{d}} 
													     \frac{e^{-(\sum_{j=1}^{p}u^{2}_{j})/4}}{\prod_{j=1}^{p}|u_{j} + x'|^{d-2}} du_{1} ... du_{p},
\end{align*}
where we use the change of variable $x'=nx$ and in the last step we use of course the fact that $c_{n} \sim n^{d-2}$ and $ V^{(n)}_{0} \longrightarrow \R^{d}$, as $n\rightarrow \infty$.

Let us set:
\[ E_{n,i}=  d_{n} D_{n,i}=  \int_{V^{(n)}_{i}}  :(\partial \psi_{n}(\frac{x}{n}))^{2}:dx = n^{d} \int_{V^{(n)}_{i}/n}  :(\partial \psi_{n}(x)^{2}:dx \]
Then we have 
\begin{align*}
  \langle E_{n,i}^{2} \rangle &= \frac{ 1}{c_{n}^{2}} n^{2d} 2!\int_{V^{(n)}_{i}/n} \int_{V^{(n)}_{i}/n} \langle \partial \varphi_{n}(x)\partial \varphi_{n}(y )\rangle^{2} dx dy \\
	                &= \frac{ 1}{c_{n}^{2}} n^{2d} |\frac{V^{(n)}_{0}}{n}|\int_{V^{(n)}_{0}/n}  \sum_{j=1}^{d}  ( \frac{\partial^{2} c_{n}}{\partial x_{j}^{2}}(x))^{2} dx
\end{align*}
To estimate the last integrals we use (\ref{SecondDerivative2}) to get:
\begin{align*}
\nonumber
\int_{V^{(n)}_{0}/n} (\frac{\partial^{2} c_{n}}{\partial x_{j}^{2}}(x))^{2} dx 
                        & \sim  \int_{V^{(n)}_{0}/n} K n^{2d} \frac{dx}{2(2\pi)^{d}} \int_{\R^{d}}\int_{\R^{d}} \frac{2(u_{j}^{2}+1) (v_{j}^{2}+1) 
												                              e^{-(u^{2}+v^{2})/4} }{|u+nx|^{d-2} |v+nx|^{d-2}} du dv \\
											  & \sim     K n^{2d} n^{-d} \int_{V^{(n)}_{0}} \frac{dx'}{2(2\pi)^{d}} \int_{\R^{d}}\int_{\R^{d}} \frac{2(u_{j}^{2}+1) (v_{j}^{2}+1) 
												                              e^{-(u^{2}+v^{2})/4} }{|u+x'|^{d-2} |v+x'|^{d-2}} du dv \\
											    & \sim     K' n^{d}
\end{align*}
\noindent
where we set $x'=nx$ and the convergence of the intergrals w.r.t. $x'$ can be seen by the change of variables $u'=u+x, v'=v+x$. Notice that we have also:
\begin{equation}
\label{IntegralSecondDeriv}
\int_{V^{(n)}_{i}} (\frac{\partial^{2} c_{n}}{\partial x_{j}^{2}}(x))^{2} dx \sim  K n^{d}\int_{n V^{(n)}_{0}} \frac{dx'}{2(2\pi)^{d}} \int_{\R^{d}}\int_{\R^{d}} \frac{2(u_{j}^{2}+1) (v_{j}^{2}+1)  e^{-(u^{2}+v^{2})/4} }{|u+x'|^{d-2} |v+x'|^{d-2}} du \sim K' n^{d},
\end{equation}
the integral is taken this time on $n V^{(n)}_{0}$.

With this and the fact that $c_{n}\sim Kn^{d-2}$ and $|V^{(n)}_{i}/n|=n^{-d}$ , we obtain:
\begin{equation}
  \langle E_{n,i}^{2} \rangle \sim  K \frac{ 1}{c_{n}^{2}} n^{2d} n^{-d}  n^{d}  \sim  K  \frac{ 1}{n^{2(d-2)}} n^{2d} \sim K n^{4}
\end{equation}
\noindent
With $d_{n}=n^{2}$ we see that:
\[ \langle D_{n,i}^{2} \rangle = \frac{1}{d_{n}^{2}} \langle E_{n,i}^{2} \rangle \]
converges to some $K>0$ as $n \longrightarrow \infty$.

The point (3) of the proposition can be proved by similar calculations: this time we deal with expressions like: $ \langle:\psi_{n}^{4}(x_{1}): ...:\psi_{n}^{4}(x_{p}):\rangle$ and
in the integrals we will get give rise to terms like  $\langle\psi_{n}(x_{1})\psi_{n}(x_{2})\rangle^{4} ... \langle\psi_{n}(x_{p-1})\psi_{n}(x_{p})\rangle^{4}$, the integrals of each couple
can be made independently of the others and this reduces the calculations to the case $p=2$.
Finally, the point (4) is easy, we omit the details. $\Box$

\noindent
Next, we have a kind of law of large numbers (LLN) for the continuous field $\varphi_{n} (x)$ in a finite volume, which implies LLNs for the arrays $ I_{n,i}, M_{n,i}, D_{n,i}, X_{n,i}$

\begin{prop}\label{PropLLN}
We have the following law of large numbers of the continuous field $\psi_{n}$:
\begin{align*}
\nonumber
 \lim_{n\longrightarrow \infty} \int_{V} :\psi_{n}^{4}(x): dx &= \lim_{n\longrightarrow \infty} \frac{1}{n^{d}} \sum_{i=1}^{n^{d}} I_{n,i} = 0, \; \mu_{0}  \; a.e  \\
 \lim_{n\longrightarrow \infty} \int_{V} :\psi_{n}^{2}(x): dx &= \lim_{n\longrightarrow \infty} \frac{1}{n^{d}} \sum_{i=1}^{n^{d}} M_{n,i} = 0, \; \mu_{0}  \; a.e  \\
 \lim_{n\longrightarrow \infty} \int_{V} \frac{1}{d_{n}}:(\partial \psi_{n})^{2}(x): dx &=  \lim_{n\longrightarrow \infty} \frac{1}{n^{d}}\sum_{i=1}^{n^{d}} D_{n,i} = 0, \; \mu_{0} \; a.e
\end{align*}
\end{prop}
{\bf Proof.}
We have:
\begin{align*}
\nonumber
   \langle (\int_{V} :\psi_{n}^{4}(x): dx)^{2}\rangle = \langle I_{n}^{2}\rangle & = \frac{1}{c_{n}^{4}}  \langle (\int_{V}:\varphi_{n}^{4}(x): dx)^{2} \rangle  \\
                                                          &= \frac{4!}{c_{n}^{4}}  \langle \int_{V}\int_{V} \langle \varphi_{n}(x)\varphi_{n}(y) \rangle^{4} dx dy  \\
		                                               				&= V \frac{4!}{c_{n}^{4}}  \int_{V} c_{n}(x)^{4} dx \\
			                                               			&\sim  K V 4! \frac{n^{3d-8}}{(n^{d-2})^{4}}   \\
				                                               		&\sim  K V 4! \frac{1}{n^{d}}
\end{align*}
which implies that $\sum_{n}  \langle I_{n}^{2}\rangle < \infty$ and that $I_{n}\longrightarrow 0, \; a.e$ (by the Borel-Cantelli lemma).
In the same way we have 
\begin{eqnarray*}
\nonumber
 \langle M_{n}^{2}\rangle &=& V \frac{1}{c_{n}^{2}}  \int_{V} c_{n}(x)^{2} dx \\
						&\sim & K V 4! \frac{\log n }{(n^{d-2})^{2}} \; {\rm if} \; d=4 \; {\rm and} \; \;  K V 4! \frac{n^{d-4} }{(n^{d-2})^{2}} \; {\rm if} \; d\geq 5   \\
\end{eqnarray*}
which implies that $\sum_{n}  \langle M_{n}^{2}\rangle < \infty$ and that $M_{n}\longrightarrow 0, \; a.e$
\\
As for $D_{n}$ we have 
\begin{align*}
\nonumber
 \langle D_{n}^{2}\rangle & = \frac{1}{d_{n}^{2}} \langle (\int_{V}  :(\partial \psi_{n})^{2}(x): dx)^{2}\rangle \\
              &= V \frac{2!}{c_{n}^{2} d_{n}^{2}} \int_{V} \sum_{j=1}^{d}  ( \frac{\partial^{2} c_{n}}{\partial x_{j}^{2}}(x))^{2} dx \\
							    &= V \frac{2!}{c_{n}^{2} d_{n}^{2}} d K n^{d} \\
									&= V \frac{2!}{n^{2(d-2)} n^{2}} d K n^{d} \\
									&= \frac{K'}{n^{d-2}}
\end{align*}
where we have used (\ref{IntegralSecondDeriv}). 
\noindent
And as before, $\sum_{n}   \langle (D_{n})^{2}\rangle < \infty$ for $d\geq 4$ which implies the almost sure convergence $D_{n}\longrightarrow 0, \; a.e.$ $\Box$


\subsection{A lower bound of large deviations probabilities for dependent sequences and arrays}

Our motivation is to obtain a lower bound for the sequence ${\cal A}_{n}^{(1)}$ like:
$$ \lim\inf_{n\longrightarrow\infty} \frac{1}{n}\log P^{(n)}( {\cal A}_{n}^{(1)}= \frac{1}{n}\sum_{i=1}^{n}X_{n,i} \in C)\geq - \inf_{{\rm int} C} I(x)$$

\noindent
Usual results in large deviations can not be applied in our case ; they often deal with the i.i.d case or dependent sequence (whose mean is $L_{n}$) with specific conditions. Few results address the case of arrays of rv's.
The Gartner-Ellis theorem can not be used: in this theorem, the LDP is deduced from conditions on the limit:
\begin{equation}
\label{GE}
 \lim_{n\longrightarrow \infty} \frac{1}{n}\log E e^{\lambda {\cal A}_{n}^{(1)}} = J (\lambda)
\end{equation}
which is supposed to exist with some properties. But this is precisely what we are looking for. In fact, we take the opposite direction, by seeking a LDP to be satisfied by ${\cal A}_{n}^{(1)}$ we wish to get an estimate of the limit  (\ref{GE}).

Lower bounds of large deviations probabilities stated in different or more general forms than those currently used
in large deviation theory are useful for many applications. One of these forms can be found in Bahadur, Zabell and Gupta \cite{BZG}
which contains some interesting examples ; we shall use a formulation given in \cite{AA1}, \cite{AA2} which deals with the i.i.d random variables case ; for the sake of clarity, we reproduce it here with its short proof. 
This formulation will be generalized to arrays of dependent random variables (Proposition \ref{PropLB2}), and the i.i.d proof will be adapted to that purpose. But to get a utilizable lower bound for the proof of Theorem \ref{main-th}, further intermediate results will be needed.

\begin{prop}\label{PropLB1}
Let $B$ be a Banach space, ${\cal B}$ a $\sigma$-field on $B$, $\Omega= B^{\N}, \; {\cal A}= {\cal B}^{\otimes \N}$ and $X^{i}: \Omega\longrightarrow B$ the
coordinate maps. Suppose that we have also the following data:\\
$\bullet$ A probability measure $P$ on $(B,{\cal B})$\\ 
$\bullet$ A probability measure $Q$ on $(B,{\cal B})$ and $C\subset B$ be such that $A_{n}:= \{\frac{1}{n} \sum_{1}^{n} X^{i} \in C \} $ is ${\cal A}$-measurable and
\begin{equation}
\label{propLB1}
 \lim_{n\longrightarrow \infty} Q^{\otimes n}(\frac{1}{n}\sum_{i=1}^{n}
X^{i} \in C) =1 .
\end{equation}
Then for all probability measures $P\in M_{1}(B)$ we have
\begin{equation}
\label{LB4}
\lim\inf_{n\longrightarrow\infty} \frac{1}{n}\log P^{\otimes n}( 
\frac{1}{n}\sum_{i=1}^{n}X^{i} \in C)\geq - D(Q|| P ) .
\end{equation}
\end{prop}

\noindent
{\bf Proof.}\\
We recall here the proof given in \cite{AA2}: we may  assume that $ Q\ll P$ (otherwise (\ref{LB4}) is obvious) and
observe that 
\begin{equation}
\label{pf11}
\int 1_{A} dP \geq \int 1_{A} \frac{1}{dQ/dP} dQ, \; \; A\in {\cal B}.
\end{equation}
(\ref{pf11}) is clearly verified if $P\ll Q$; otherwise, taking $Q_{\alpha}=\alpha Q +(1-\alpha) P$, we have
\[ P(A) \geq \alpha \int 1_{A} \frac{1}{dQ_{\alpha}/{dP}} dQ,  \]
and we get (\ref{pf11}) by using  Fatou's lemma.
Now, applying (\ref{pf11})  to $ P^{\otimes n}$ and $ Q^{\otimes n}$ and using
the Jensen inequality, we have
\begin{align*}
\log  P^{\otimes n}(A_{n}) &\geq \log \int 1_{A_{n}} 
                 \frac{1}{ d  Q^{\otimes n}/dP^{\otimes n}} d Q^{\otimes n} \\
        &\geq  -\frac{1}{Q^{\otimes n}(A_{n})}\int 1_{A_{n}} \log\frac{ d  Q^{\otimes n}}{
 dP^{\otimes n}} d Q^{\otimes n}\\
        &= -n \frac{1}{Q^{\otimes n}(A_{n})}\int 1_{A_{n}} \log \frac{d Q}{d P} 
  d Q^{\otimes n}.
\end{align*}
Consequently
\begin{equation}
 \frac{1}{n}\log P^{\otimes n}(A_{n} )\geq -\frac{1}{Q^{\otimes n}(A_{n})}
    \int 1_{A_{n}} \log\frac{d Q}{d P}
d Q^{\otimes \N};
\end{equation}
hence, using (\ref{propLB1}) and Lebesgue's theorem we get (\ref{LB4}).
 $\Box$

\noindent
The next proposition is a generalization of this result to the case of array of dependent random variables. 
Let $B$ be a Banach space and ${\cal B}$ a $\sigma$-field on $B$ (e.g. the Borel  $\sigma$-field). 

\begin{prop}\label{PropLB2}
Let $B$ be a Banach space and $X_{N,i}, i=1, ..., N $ an array of $B$-valued rv's. We suppose that the $X_{N,i}$ have the same law $P_{N}$ but are not necessarily independent. Let $P_{N}^{(N)}$ be the law of  $(X_{N,i}, i=1, ..., N) $, defined as a measure on $B^{N}$. We denote by ${\cal A}$ the $\sigma$-field $\sigma(X_{N,i}, i=1, ..., N, N\geq 1)$. Suppose that we have also the following data:

$\bullet$ The probability measures $P_{N}^{\otimes N}$ on $(B^{N},{\cal B}^{N})$ with the same marginals $P_{N}$, so that the $(X_{N,i}, i=1, ..., N) $ are i.i.d under $P_{N}^{\otimes N}$. 

$\bullet$ A probability measure $Q_{N}$ on $(B,{\cal B})$ ; we also consider the the pm's $Q_{N}^{\otimes N}$ on $(B^{N},{\cal B}^{N})$, under which the $(X_{N,i}, i=1, ..., N) $ are i.i.d  with the law $Q_{N}$.

$\bullet$ A subset $C\subset B$ be such that $A_{N}:= \{\frac{1}{N} \sum_{1}^{N} X_{N,i} \in C \} $ is ${\cal A}$-measurable and

\begin{equation}
\label{prop21}
 \lim_{n\longrightarrow \infty} Q^{\otimes N}_{N}(   \frac{1}{N}  \sum_{i=1}^{N} X_{N,i}\in C) =1 .
\end{equation}
If:
\begin{equation}
\label{prop211}
 \lim_{n\longrightarrow \infty} D( Q_{N} || P_{N} ) = D \; { \rm exists},
\end{equation}
then we have
\begin{equation}
\label{prop22}
\lim\inf_{n\longrightarrow\infty} \frac{1}{N}\log P^{(N)}_{N}(  \frac{1}{N}  \sum_{i=1}^{N} X_{N,i} \in C)\geq - D.
\end{equation}
\end{prop}
For the proof of this proposition, we need some intermediate results stated in the following lemma and propositions:
\begin{lm}\label{lmLB1}
Let $ P^{\otimes N}=\otimes_{i=1}^{N} P_{i}$ be a product measure on a space $E= E_{1} \times ... \times E_{N}$ and $ Q^{N}$ a probability measure on $E$ such that each marginal
$Q_{i}= Q^{N} |_{E_{i}}$ is absolutely continuous w.r.t the corresponding marginal $P_{i}$ of $P$  ( $Q_{i} \ll P_{i}$). We also suppose that $ Q^{N}$ is supported by the open sets of $E$. Then $ Q^{N} \ll P^{N} $ 
\end{lm}
{\bf Proof.}\\
Let $ A = A_{1} \times ... \times A_{N}$ be such that $P^{N}(A)=0$. Then there exists some $j$ such that $P_{j} (A_{j})=0 $. On the other hand, 
\[ Q^{N}(A) \leq Q^{N}(E_{1} \times ... \times A_{j} \times ... \times E_{N}) = Q_{j}(A_{j}).  \]
\noindent
And since $Q_{j} \ll P_{j}$ and $P_{j} (A_{j})=0 $ we have $Q_{j}(A_{j})=0$ and $Q^{N}(A)=0$.
\\
\noindent
Now, let $A \subset E$ be an open set such that $P^{N}(A)=0$. Then for all rectangular set $ A' = A_{1} \times ... \times A_{N} \subset A $, we have  $P^{N}(A')=0$ which, by the previous argument, implies that $Q^{N}(A')=0$ and as $A$ is a countable union of rectangles, this gives  $Q^{N}(A)=0$.
 $\Box$
\\
\\
\noindent
The following proposition is inspired from Bahadur and Raghavachari \cite{BR} Theorem 1 p. 133, which provides an interesting property of the limit 
of $(1/n)\log dQ/dP|_{{\cal B}^{n}} $ with ${\cal B}^{n}$ a sequence of $\sigma-$fields:

\begin{prop}\label{PropLB3}
Let $P,Q$ be two probability measures on a measurable space $(B,{\cal B})$ and ${\cal B}^{n} \subset {\cal B}$ a sequence of $\sigma-$fields.
\\
Suppose that $Q$ is absolutely continuous w.r.t $P$ on ${\cal B}^{n}$, with 
\begin{equation}
\label{BR1}
R_{n}= \frac{dQ}{dP}\arrowvert_{{\cal B}^{n}}
\end{equation}
\noindent
Then:
\begin{equation}
\label{BR2}
K:= \lim\inf_{n\longrightarrow\infty} \frac{1}{n}\log R_{n} \geq 0  \; Q \; a.e.
\end{equation}
\end{prop}

\noindent
{\bf Proof.}\\
Let $K_{n}=(\log R_{n})/n$. Then, for $\varepsilon > 0$, $K_{n} < -\varepsilon$ iif =$ R_{n} < \exp (-n \varepsilon) $ and:
\begin{align*}
\nonumber
Q (K_{n} < -\varepsilon) &= \int_{K_{n} < -\varepsilon} R_{n} dP \\
                      &= \int_{R_{n}< \exp (-n \varepsilon)} R_{n} dP \\
											&\leq \int e^{-n \varepsilon} dP \\
											&= e^{-n \varepsilon}
\end{align*}
Now let $S_{n}$ be the event $ \{K_{n} < -\varepsilon \}= \{ -K_{n}  > \varepsilon \}$. From the last inequality we get $ \sum_{n\geq 1} Q (S_{n}) < \infty$, which implies that
$ Q(\limsup S_{n}) = 0$. And:
\begin{align*}
\nonumber
 \limsup S_{n} &= \{ \limsup -K_{n}  > \varepsilon \} \\
               &= \{ -\liminf K_{n}  > \varepsilon \}  \\
               &= \{ \liminf K_{n}  < - \varepsilon \}
\end{align*}
This shows that $ Q (\{ \liminf K_{n}  < -\varepsilon \})=0, \forall \varepsilon > 0$  and $ Q (\{ \liminf K_{n}  \geq 0 \})=1$. $\Box$

We associate to the array of random variables $X_{N,i}, i=1, ..., N $, the sequence of $\sigma-$fields ${\cal B}^{N} = \sigma (X_{N,i}, i=1, ..., N)$.
The sequence of probability measure $P^{(N)}_{N}$ defined on ${\cal B}^{n}$ can be extended to a probability measures $P^{\infty}$ on 
${\cal B} = \sigma (X_{N,i}, i=1, ..., N, N \geq 1)$ such that $P^{\infty}|_{{\cal B}^{N}}=P^{(N)}_{N} $.

\begin{prop}\label{PropLB4}
With the notations of Proposition \ref{PropLB2}, let:
\begin{equation}
\label{prop331}
R_{n}= \frac{d P^{(N)}_{N}}{dP^{\otimes N}_{N}}
\end{equation}
Then:
\begin{equation}
\label{prop32}
K:= \lim\inf_{n\longrightarrow\infty} \frac{1}{n}\log R_{n} \geq 0 \; \; \; P^{\infty} \;  a.e.
\end{equation}
\end{prop}
{\bf Proof.}\\
In the same way, the sequence $Q^{\otimes N }_{N}$ defined on ${\cal B}^{N}$ can be extended to a probability measure $Q^{\infty}$ on 
${\cal B} = \sigma (X_{N,i}, i=1, ..., N, N \geq 1)$ such that $Q^{\infty}|_{{\cal B}^{N}}=Q^{\otimes N }_{N} $. 
\\
On  ${\cal B}^{N} $ we have $P^{\infty}|_{{\cal B}^{N}} \ll P^{\otimes N}_{N}$. The previous proposition can be applied and we have: 

\[ K:= \lim\inf_{n\longrightarrow\infty} \frac{1}{N}\log R_{N} \geq 0  \; \; P^{\infty} \;  a.e. \; \; \Box \]

\noindent
{\bf Proof of Proposition \ref{PropLB2}}\\
By the previous lemma \ref{lmLB1}, we have: $P_{N}^{(N)} \ll P_{N}^{\otimes N}$. Then, by the same arguments as those of the beginning of the proof of Proposition \ref{PropLB1}:
\[ P_{N}^{(N)} ( A_{N}) \geq \int 1_{A_{N}} dP \geq \int 1_{A_{N}} \frac{d P_{N}^{(N)}}{d P_{N}^{\otimes N}} \frac{1}{dQ_{N}^{\otimes N} /d P_{N}^{\otimes N}}dQ_{N}^{\otimes N} \]
By the Jensen's inequality, we have 
\begin{align*}
\nonumber
\log P_{N}^{(N)} ( A_{N}) &\geq  \int 1_{A_{N}} \log \frac{d P_{N}^{(N)}}{d P_{N}^{\otimes N}} dQ_{N}^{\otimes N}  - 
                                  \int 1_{A_{N}} \log \frac {d Q_{N}^{\otimes N}}{dP_{N}^{\otimes N}} dQ_{N}^{\otimes N}\\
										      &=  \int 1_{A_{N}} \log \frac{d P_{N}^{(N)}}{d P_{N}^{\otimes N}} dQ_{N}^{\otimes N}  - 
													     \int \log \frac {d Q_{N}^{\otimes N}}{dP_{N}^{\otimes N}} dQ_{N}^{\otimes N} +\\
													 &+ 		\int 1_{A_{N}^{c}} \log \frac {d Q_{N}^{\otimes N}}{dP_{N}^{\otimes N}} dQ_{N}^{\otimes N}\\
													&= T_{N}^{(1)} + T_{N}^{(2)} + T_{N}^{(3)}
\end{align*}
We have to estimate each of these 3 terms. For the second one we have:
\begin{eqnarray*}
\nonumber
\frac{1}{N} T_{N}^{(2)} &=& \frac{1}{N} N \int \log \frac {d Q_{N}}{dP_{N}} dQ_{N}\\
                       &=&  D( Q_{N} || P_{N} )\\
											 &\longrightarrow & D
\end{eqnarray*}
as $N \longrightarrow \infty$ by the assumption of the proposition. As for the 3d term, using $ x\log x \geq x -1$:
\begin{align*}
\frac{1}{N} T_{N}^{(3)} &= \frac{1}{N} \int 1_{A_{N}^{c}} \log \frac {d Q_{N}^{\otimes N}}{dP_{N}^{\otimes N}} \frac{dQ_{N}^{\otimes N}}{dP_{N}^{\otimes N}}dP_{N}^{\otimes N} \\
                       &\geq  \frac{1}{N} [ \int 1_{A_{N}^{c}} \frac {d Q_{N}^{\otimes N}}{dP_{N}^{\otimes N}} dP_{N}^{\otimes N} - \int 1_{A_{N}^{c}} dP_{N}^{\otimes N} ]   \\
											 &=  \frac{1}{N} ( Q_{N} (A_{N}^{c}) - P_{N} (A_{N}^{c}) )\\
											&\longrightarrow  0
\end{align*}
as $ N \longrightarrow \infty$.
\noindent
Next we turn to the proof of:
\[  T^{(1)}=\liminf_{N \longrightarrow \infty} \frac{1}{N} T_{N}^{(1)} = 0  \]

We have:
\[  \liminf_{N \longrightarrow \infty} \frac{1}{N} T_{N}^{(1)} = \liminf_{N \longrightarrow \infty} \int 1_{A_{N}} \frac{1}{N} \log \frac{d P_{N}^{(N)}}{d P_{N}^{\otimes N}} dQ_{N}^{\otimes N}  \]

\noindent
By proposition \ref{PropLB4}, we have $ \liminf_{N \longrightarrow \infty} (1/N)\log (d P_{N}^{(N)}/d P_{N}^{\otimes N}) \geq 0 \; \; P^{(\infty)} \; a.e.$ ; but it does not imply directly that $T^{(1)}=0$, because the integration in the left hand side is w.r.t $dQ_{N}^{\otimes N}$ and we have not proved that it is absolutely continuous w.r.t $P^{(\infty)}$. We can proceed as follows: By Proposition \ref{PropLB4} there exists a subset $B_{1}$ with $P^{(\infty)}(B_{1}) = 1$ such that:

\[ \liminf_{N \longrightarrow \infty} \frac{1}{N} \log \frac{d P_{N}^{(N)}}{d P_{N}^{\otimes N})} (\omega) \geq 0 \; \; \forall \omega \in B_{1} \]
\noindent
In the integration performed in the begining of the proof we use $B_{1}$ instead of integrating in the whole space:
\begin{align*}
\nonumber
\log P_{N}^{(N)} ( A_{N}) &\geq  \frac{1}{ Q_{N}^{\otimes N}(A_{N})} [\int 1_{B_{1}} 1_{A_{N}} \log \frac{d P_{N}^{(N)}}{d P_{N}^{\otimes N}} dQ_{N}^{\otimes N}  - 
                                  \int 1_{B_{1}} 1_{A_{N}} \log \frac {d Q_{N}^{\otimes N}}{dP_{N}^{\otimes N}} dQ_{N}^{\otimes N} ]\\
                           &\geq   \frac{1}{ Q_{N}^{\otimes N}(A_{N})} [ \int 1_{B_{1}} 1_{A_{N}} \log \frac{d P_{N}^{(N)}}{d P_{N}^{\otimes N}} dQ_{N}^{\otimes N} -
                                      \int \log \frac {d Q_{N}^{\otimes N}}{dP_{N}^{\otimes N}} dQ_{N}^{\otimes N} \\
                            &+  \int 1_{(A_{N} \cap B_{1})^{c}} \log \frac {d Q_{N}^{\otimes N}}{dP_{N}^{\otimes N}} dQ_{N}^{\otimes N} ]\\
														&= T_{N}^{'(1)} + T_{N}^{(2)} + T_{N}^{'(3)}
\end{align*}
Then
\begin{equation}
\label{LBineq}
\log P_{N}^{(N)} ( A_{N}) \geq - N D( Q_{N} || P_{N}) +  T_{N}^{'(1)} + T_{N}^{'(3)}
\end{equation}

\noindent
By Proposition \ref{PropLB4}:
$$\liminf_{N \longrightarrow \infty} 1_{B_{1}} (1/N)\log (d P_{N}^{(N)}/d P_{N}^{\otimes N}) \geq 0$$
hence $\liminf T_{N}^{'(1)}/N \geq 0$. As for $T_{N}^{'(3)}$, using again $ x\log x \geq x -1$:
\begin{align*}
\nonumber
\frac{1}{N} T_{N}^{'(3)} &= \frac{1}{N} \frac{1}{ Q_{N}^{\otimes N}(A_{N})} \int 1_{(A_{N} \cap B_{1})^{c}} \log \frac {d Q_{N}^{\otimes N}}{dP_{N}^{\otimes N}} 
                             \frac{dQ_{N}^{\otimes N}}{dP_{N}^{\otimes N}}dP_{N}^{\otimes N} \\
                       &\geq  \frac{1}{N} \frac{1}{ Q_{N}^{\otimes N}(A_{N})} [ \int 1_{(A_{N} \cap B_{1})^{c}} \frac {d Q_{N}^{\otimes N}}{dP_{N}^{\otimes N}} 
											              dP_{N}^{\otimes N} - \int 1_{(A_{N} \cap B_{1})^{c}} dP_{N}^{\otimes N} ]   \\
											 &=  \frac{1}{N} \frac{1}{ Q_{N}^{\otimes N}(A_{N})} ( Q_{N} (A_{N} \cap B_{1})^{c}) - P_{N} (A_{N} \cap B_{1}^{c})\\
											&\longrightarrow  0
\end{align*}
as $ n \longrightarrow \infty$.
\noindent
Finally, as before, $T_{N}^{(2)}/N \longrightarrow D$. $\Box$
\\
The next proposition is the initial motivation of the previous general lower bound of large deviations.

\begin{prop}
\label{PropLB5}
Let $\lambda_{n} \geq 0, \alpha_{n}, \beta_{n}$ be a sequence of reals that converge to limits $\lambda, \alpha, \beta$ that are finite and at least one of them is not null.
We suppose that $\lambda \geq 0$. Then, with
\begin{equation}
X_{n,i} = \lambda_{n} \int_{V^{(n)}_{i}} :\psi_{n}^{4}(\frac{x}{n}): dx - \alpha_{n}\int_{V^{(n)}_{i}} :\psi_{n}^{2}(\frac{x}{n}): dx - \beta_{n}  \frac {1}{d_{n}} \int_{V^{(n)}_{i}}  :(\partial \psi_{n})^{2}(\frac{x}{n}): dx
\end{equation}
there exists a function $I(h)$ such that:

\begin{equation}
\label{LB5}
\lim\inf_{n\longrightarrow\infty} \frac{1}{n}\log P^{(N)}_{N}(  \frac{1}{N}  \sum_{1}^{N} X_{N,i} \in C)\geq -\inf_{h \in {\rm int} C} I(h)
\end{equation}
In other words:
\begin{equation}
\label{LB6}
\lim\inf_{n\longrightarrow\infty} \frac{1}{n}\log P^{(N)}_{N}(L_{n}(\psi) \in C)\geq -\inf_{h \in {\rm int} C} I(h)
\end{equation}
With 
\[ L_{n}(\psi) = \lambda_{n} \int_{V} :\psi_{n}^{4}(x): dx + \alpha_{n}\int_{V} :\psi_{n}^{2}(x): dx + \beta_{n}\frac{1}{d_{n}}\int_{V}:(\partial \psi_{n})^{2}(x): dx \]
The function $I$ has the following properties:
\\
$\bullet$ $I(0)=0$ and the definition domain of $I$ includes the interval $[0,b], b>0 $ and for a given constant $K>0$, $b$ may be chosen such that for all $h\in [0, b]: I(h) \leq K h$.
\end{prop}

\noindent
The proof of this proposition requires some intermediate results stated in the following lemma and propositions. 
\\
\\
{\bf Notations.} For 
\[ X_{n,i}= \lambda_{n}I_{n,i} - \alpha_{n} M_{n,i} - \beta_{n} D_{n,i}  
\]
let us set:
\[ L_{n}^{I}(\theta) = E e^{-\lambda_{n} \theta I_{n,i}},  \; \; \Lambda_{n}^{I}(\theta)= \log E e^{-\lambda_{n} \theta I_{n,i}} \]
\[ L_{n}^{M}(\theta) = E e^{-\theta \alpha_{n} M_{n,i}}, \; \;  \Lambda_{n}^{M}(\theta)= \log E e^{-\theta \alpha_{n} M_{n,i}} \]
\[ L_{n}^{D}(\theta) = E e^{-\theta \beta_{n} D_{n,i}}, \; \;  \Lambda_{n}^{I}(\theta)= \log E e^{-\theta \beta_{n}D_{n,i}} \]
\[ L_{n}^{X}(\theta) = E e^{-\theta X_{n,i}}, \; \;  \Lambda_{n}^{X}(\theta)= \log E e^{-\theta X_{n,i}} \]

\begin{lm}\label{lmLB2}
Suppose that that the sequences $\lambda_{n} \geq 0, \alpha_{n}, \beta_{n}$ have limits $\lambda, \alpha, \beta$. Then:
\\
$\bullet$ The functions $L_{n}^{I}(\theta),\Lambda_{n}^{I}(\theta)$ are uniformly bounded in every closed interval of $\R^{+}$. The same result holds for $L_{n}^{M}(\theta), \Lambda_{n}^{M}(\theta)$ and $L_{n}^{D}(\theta), \Lambda_{n}^{D}(\theta)$ in in any closed subinterval of $[0, 1/2\alpha[ $ and $[0, 1/2\beta[ $ respectively.
\\
$\bullet$ The derivatives of $L_{n}^{I}(\theta), \Lambda_{n}^{I}(\theta)$ are uniformly bounded in every closed interval of $\R^{+}$. Those of $L_{n}^{M}(\theta), \Lambda_{n}^{M}(\theta)$ are uniformly bounded in every closed interval of  $[0, 1/(2\alpha)[ $.
\\
$\bullet$ The derivatives of $ L_{n}^{D}(\theta), \Lambda_{n}(\theta)$ are uniformly bounded in every closed interval of $[0, 1/(2\beta)[ $.
\\
$\bullet$ The functions $L_{n}^{X}(\theta), \Lambda_{n}^{X}(\theta)$ and their first and second derivatives are bounded in any closed subinterval of $[0, 1/4(\max(\alpha,\beta))[$.
\end{lm}
\noindent
{\bf Proof.}\\
To begin with $L_{n}^{I}$, let us remark that:
\begin{align*}
\nonumber
E e^{-\lambda_{n} \theta I_{n,i}} &= \int_{0}^{+\infty} Pr(e^{-\lambda_{n} \theta I_{n,i}} > t) dt\\
                     &\leq  1+\int_{1}^{+\infty} Pr(e^{-\lambda_{n} \theta I_{n,i}} > t) dt \\
										&\leq  1+\int_{1}^{+\infty} Pr(-I_{n,i} > \frac{\log t}{\lambda_{n}\theta}) dt \\
										&\leq  1+\int_{1}^{+\infty} Pr(-\int_{V} :\psi_{n}^{4}(x): dx  > \frac{\log t}{\lambda_{n}\theta}) dt \\
										&\leq  1+\int_{1}^{+\infty} Pr(-:\psi_{n}^{4}(x):  > \frac{\log t}{\lambda_{n}\theta} \; \; {\rm for \; some } \; x ) dt
\end{align*}
The last inequality is due to the fact that $|V|=1$. On the other hand, we have always $-:\psi_{n}^{4}(x): \leq 3$, then 
\begin{eqnarray*}
E e^{-\theta I_{n,i}}&\leq & 1+\int_{1}^{e^{3 \lambda_{n} \theta}} Pr(-:\psi_{n}^{4}(x):  > \frac{\log t}{\theta} \; \; {\rm for \; some } \; x  ) dt\\
                     &\leq & 1+e^{3 \lambda_{n} \theta}
\end{eqnarray*}
\noindent
so that $L_{n}^{I}(\theta)$ is bounded in any closed interval of $\R^{+}$, because of the convergence of $\lambda_{n}$ . As for the the derivatives: 

\[ {L_{n}^{I}}'(\theta)= E \lambda_{n} I_{n,i} e^{-\lambda_{n}\theta I_{n,i}}, \; \;  L_{n}^{I (2)}(\theta)= \lambda_{n}^{2}E  I_{n,i}^{2} e^{-\lambda_{n}\theta I_{n,i}}, \] 

\noindent
their boundedness follows from that of  $L_{n}^{I}(\theta)$ and $E (I_{n,i})^{4}$ by proposition \ref{PropArrays} and the Cauchy-Schwarz inequality. 
\\
The boundedness of $\Lambda_{n}(\theta) = \log L_{n}^{I}(\theta) $ follows from that of $L_{n}^{I}(\theta)$ and the fact that $L_{n}^{I}(\theta) \geq \exp[-\theta E X_{n,i}] = 1$ (Jensen's inequality).
\\
Let us turn to the case of $L_{n}^{M}(\theta)= E \exp(-\alpha_{n}\theta M_{n,i})$. To simplify the notations, we suppose first that $\alpha_{n}\equiv 1$. We have: 
\begin{align*}
\nonumber
E e^{-\theta M_{n,i}} &= \int_{0}^{+\infty} Pr(e^{-\theta M_{n,i}} > t) dt \\
                     &\leq  e+\int_{e}^{+\infty} Pr(e^{-\theta M_{n,i}} > t) dt \\
										&\leq  e+\int_{e}^{+\infty} Pr(\int_{V} :\psi_{n}^{2}(x): dx  > \frac{\log t}{\theta}) dt \\
										&\leq  e+\int_{e}^{+\infty} Pr(V :\psi_{n}^{2}(x):  > \frac{\log t}{\theta} \; \; {\rm for \;  some } \; x ) dt \\
										&\leq  e+\int_{e}^{+\infty} Pr( \psi_{n}^{2}(x)  > \frac{\log t}{V \theta} \; \; {\rm for \; some } \; x ) dt \\
										&\leq  e+\int_{e}^{+\infty} Pr( |N|  > (\frac{\log t}{V \theta})^{1/2} ) dt
\end{align*}
where $N$ is a normal gaussian random variable. The tail probabilities of $N$ have bounds like:
\[ Pr( N  > y) = \frac{1}{\sqrt{2\pi}} \int_{y}^{+\infty} e^{-x^{2}/2}dx \leq \frac{1}{\sqrt{2\pi}} \int_{y}^{+\infty} \frac{x}{y} e^{-x^{2}/2}dx=\frac{e^{-y^{2}/2}}{y\sqrt{2\pi} }. \]
We get:
\begin{eqnarray*}
E e^{-\theta M_{n,i}} &\leq&  e+ \frac{2}{\sqrt{2\pi}} \int_{e}^{+\infty} \frac{\exp[-(\log t^{1/\theta})/2]}{ (\log t^{1/\theta})^{1/2} }   \\
                     &\leq&  e+ \frac{2}{\sqrt{2\pi}} \int_{e}^{+\infty} \frac{ dt}{ t^{1/(2\theta)}(\log t^{1/\theta})^{1/2}  }
\end{eqnarray*}

The last integral is known to converge iif $1/(2\theta) > 1$, i.e. $\theta  < 1/2$ which proves the uniform boundedness (in $n$) of $E \exp[-\theta M_{n,i}]$ and in the same way that of $E e[-\theta D_{n,i}] $ in any closed subinterval of $[0,1/2[$. We have the same result for their derivatives with the same arguments as for $E \exp[-\theta I_{n,i}]$.
\\
Now if $\alpha_{n} \longrightarrow \alpha$, we replace the previous Gaussian variable $N$ by a Gaussian variable $G_{n}=\alpha_{n} N$ with variance $\alpha_{n}$. The factor $1/(2\theta)$ in the proof will be replaced by $1/(2\alpha_{n}\theta)$ and by a simple limit argument, the boundedness of $L_{n}^{M}(\theta)$ is valid if $1/(2\alpha\theta) > 1 $ i.e. $\theta < 1/(2\alpha)$. The case of $L_{n}^{D}$ is treated in the same way.
\\
As for the case of $L_{n}^{X}(\theta)$, we use the Holder inequality with:

$$1/p_{1}= 2\varepsilon, 1/p_{2}= 1/p_{3}= 1/2- \varepsilon, p_{2}=p_{3}= (2/(1-2\varepsilon)) {\rm \; for \; an \;} \varepsilon > 0$$
\noindent
and we have:
\begin{align*}
\nonumber
E e^{-\theta X_{n,i}} &\leq  (E e^{- p_{1} \lambda_{n}\theta I_{n,i}})^{1/p_{1}}  (E e^{ p_{2} \alpha_{n} \theta M_{n,i}})^{1/p_{2}} 
                                 (E e^{ p_{3} \beta_{n}\theta D_{n,i}})^{1/p_{3}}  \\
                      &\leq (E \exp [- (1/2\varepsilon) \lambda_{n}\theta I_{n,i}])^{2\varepsilon}  (E \exp[(2/(1-2\varepsilon)) \alpha_{n} \theta M_{n,i}])^{(2/(1-2\varepsilon))} \\
											 &  \times  (E \exp[(2/(1-2\varepsilon)) \beta_{n}\theta D_{n,i}])^{(2/(1-2\varepsilon))}
\end{align*}
In view of the first assertions of the lemma, the uniform boundedness of $L_{n}^{X}(\theta)$ is guaranteed if $ (2/(1-2\varepsilon)) \theta < \min( 1/(2\alpha), 1/(2\beta)$, that is $\theta < (1-2\varepsilon)/(\min( 1/(4\alpha), 1/(4\beta))$. As $\varepsilon$ is arbitrary we get the result for every closed interval $[0, a]$ with $ a < 1/(4 \max(\alpha, \beta))$. The same result holds of course of $\Lambda^{X}_{n}$, and the case of the derivatives of $L_{n}^{M}, L_{n}^{D}, L_{n}^{X}, \Lambda^{X}_{n}$ is proved in the same way. $\Box$
\\
\\
The following Proposition is another important piece of the proof:
\begin{prop}
\label{PropLB6}. With the assumptions of Proposition \ref{PropLB5}, let $X_{n,i}= \lambda_{n}I_{n,i} - \alpha_{n} M_{n,i} - \beta_{n} D_{n,i}$ and $\Lambda_{n}^{*}$ be the Legendre-Fenchel transform of $\Lambda_{n}$: $\Lambda_{n}^{*}(h)=\sup_{\theta \in \R}(\theta h - \Lambda_{n}(\theta))$, with $\Lambda_{n} (\theta) = E \exp[ -\theta X_{n,i}] $ . Let:
\[ \Lambda^{*}(h) = \limsup_{n\longrightarrow \infty} \Lambda_{n}^{*}(h) \]
Then, the definition domain of $\Lambda^{*}$ is not empty: there exists a real $b>0$ such that for all $h \in [0, b]$,  $ \Lambda^{*}(h)$ is finite and moreover, for a given $K>0 $, we may choose $b$ such that $ \Lambda^{*}(h) \leq K h $ for all $h \in [0, b]$. 
\end{prop}

Let us recall some facts about the Legendre-Fenchel transform $f^{*}(x)= \max_{\theta\in \R} (x\theta - f(\theta))$ in the case where $f$ is a convex function. As $\theta \longrightarrow x\theta - f(\theta)$ is concave, it will have a maximum iif $(x\theta - f(\theta))'= 0$ at some point, i.e. iif there exists a $\theta$ such that $x=f'(\theta)$. Now, suppose that $f'$ is strictly increasing and has thus an inverse at the point $x$: $\theta= f'^{-1}(x)$. Then:
\begin{eqnarray*}
\nonumber
 {f^{*}}'(x) &=&  x {f'}^{-1}(x) - f({f'}^{-1}(x)) \\
          &\Rightarrow  & {f^{*}}'(x)=  {f'}^{-1}(x) + x ({f'}^{-1}(x))'- f'(f'^{-1}(x))\times ({f'}^{-1}(x))' \\
					&= & {f'}^{-1}(x).
\end{eqnarray*}
\noindent
We shall also need the following elementary technical lemma:

\begin{lm}\label{lmLB4}
Let $Y_{n}= Y_{n}^{+} - Y_{n}^{-}$ be a sequence of real random variables splitted into their positive and negative parts. Suppose that:\\
(1) $E Y_{n} = 0$ for all $n$ \\
(2) $\liminf_{n\longrightarrow \infty} E Y_{n}^{2} = \sigma   > 0$.\\
(3) There exists $ M> 0$ such that $ E |Y_{n}|^{3} \leq M $ for all $n$. \\
Then
\begin{equation}
\label{PositivePart1}
\liminf_{n\longrightarrow \infty} E Y_{n}^{+}=\liminf_{n\longrightarrow \infty} E Y_{n}^{-} = l > 0
\end{equation}
and:
\begin{equation}
\label{PositivePart2}
\forall \varepsilon > 0 \; \exists n_{1}(\varepsilon) \ ; \; Pr ( Y_{n}^{+} > l-\varepsilon)  > 0 \; \; Pr ( Y_{n}^{-} > l-\varepsilon)  > 0 \; \; \forall n \geq n_{1}(\varepsilon).
\end{equation}
\end{lm}
\noindent
{\bf Proof.}\\
As $E X_{n} = 0$ , we have $E Y_{n}^{+} = E Y_{n}^{-}$, and $E |Y_{n}| = 2E Y_{n}^{+}$. By the Cauchy-Schwarz inequality: 

\[ E Y_{n}^{2} =  E|Y_{n}|^{1/2}|Y_{n}|^{3/2} \leq (E|Y_{n}|)^{1/2} (E Y_{n}^{3})^{1/2} \]
which yields
\[ E|Y_{n}|= 2E Y_{n}^{+} \geq \frac{E Y_{n}^{2}}{M}  \Rightarrow 2\liminf_{n\longrightarrow \infty} E Y_{n}^{+}=2\liminf_{n\longrightarrow \infty} E Y_{n}^{-} \geq \frac{\sigma}{M}. \]

This proves (\ref{PositivePart1}). As for (\ref{PositivePart2}), for a given $\varepsilon > 0$, by (\ref{PositivePart1}), there exists $n_{1}(\varepsilon)$ such that $E Y_{n}^{+} \geq l-\varepsilon/2$ for all $n \geq n_{1}(\varepsilon)$. On the other hand:

\[ E Y_{n}^{+} = E Y_{n}^{+} 1_{Y_{n}^{+}\leq l - \varepsilon} + E Y_{n}^{+} 1_{Y_{n}^{+}  > l - \varepsilon} \leq l - \varepsilon + E Y_{n}^{+} 1_{Y_{n}^{+}  > l - \varepsilon} \]

Hence for all $n \geq n_{1}(\varepsilon)$, we should have $ E Y_{n}^{+} 1_{Y_{n}^{+}  > l - \varepsilon} \geq \varepsilon/2 > 0$ and for this we must have $Pr(Y_{n}^{+}  > l - \varepsilon) > 0, \forall n \geq n_{1}(\varepsilon) $. The same result holds for $Y_{n}^{-}$. $\Box$

\begin{Rk}
The results of this lemma fail if the condition (3): $ E |Y_{n}|^{3} \leq M $ is not satisfied: consider this example: 
\begin{eqnarray*}
Y_{n} &= & x_{n}= (1-\frac{1}{n})^{-1} \times \frac{1}{\sqrt{n}} \; \; {\rm with \; probability} \;  p^{x}_{n}= 1-\frac{1}{n} \\
Y_{n} &= &x_{n}= - \sqrt{n}  \; \;   {\rm with \;  probability} \; p^{y}_{n}= \frac{1}{n}
\end{eqnarray*}
\noindent
We have $ E Y_{n}= 0$ and $ E Y_{n}^{2}= 1+ 1/n(1- 1/n^) \longrightarrow 1$. But $E Y_{n}^{+} = 1/\sqrt{n} \longrightarrow 0$. Here the condition (3) is not satisfied: $E |Y_{n}|^{3} \longrightarrow \infty$.
\\
Let us note that the estimate of positive/negative of a random variable in general is a quite complex problem.
\end{Rk}

\noindent
{\bf Proof of Proposition \ref{PropLB6}}\\
First, by lemma \ref{lmLB2}, the functions $\Lambda_{n}(\theta)$ and their first and second derivatives $\Lambda_{n}', \Lambda_{n}''$ are well defined and bounded in an interval $[0, a]$, with $a> 0$ independent of $n$. Moreover $\Lambda_{n}'$ is strictly increasing in this interval and is a bijection from $[0,a]$ to $[0, \Lambda_{n}'(a)]$. 
We can choose $a$ to be $\leq K $.
\\
Now, a crucial point is that:
\begin{equation}
\label{Defdomain1}
 l:=\liminf_{n\longrightarrow \infty} \Lambda_{n}'(a) > 0
\end{equation}
To prove (\ref{Defdomain1}), let $Y_{n}=-X_{n,i}$ and $Y_{n,i}= Y_{n,i}^{+} - Y_{n,i}^{-}$ be the positive / negative parts decomposition of $Y_{n,i}$ ; as we consider the expectations we drop the $i$ index, because, the laws of $X_{n,i}$ are independent of it. Since $xe^{x}$ is convex when $x\geq 0$ we have $ E Y_{n}^{+} e^{a Y_{n}^{+}} \geq  (E Y_{n}^{+}) e^{a EY_{n}^{+}}$ ; on the other hand $ E Y_{n}^{-} e^{- a Y_{n}^{-}} \leq E Y_{n}^{-} = E Y_{n}^{+} $, so that:
\[ E Y_{n} e^{a Y_{n}} = E Y_{n}^{+} e^{a Y_{n}^{+}} - E Y_{n}^{-} e^{- a Y_{n}^{-}} \geq E Y_{n}^{+} (e^{a EY_{n}^{+}} -1), \]

\noindent
from which we deduce that:
\[ \liminf_{n\longrightarrow \infty} E Y_{n} e^{a Y_{n}} \geq (\liminf_{n\longrightarrow \infty} E Y_{n}^{+}) (e^{a \liminf_{n\longrightarrow \infty} EY_{n}^{+}} -1), \]

\noindent
where we have used the property that $ \liminf_{n\longrightarrow \infty} f(x_{n})= f (\liminf_{n\longrightarrow \infty} x_{n}) $ provided that $f$ is continuous and increasing (this last condition is necessary), which is applied to $f(x)= x(e^{x}-1), \; x\geq 0$.
\\
We can use lemma \ref{lmLB4} for $Y_{n,i}$: its 3 conditions are satisfied by Proposition \ref{PropArrays}: $\liminf_{n\longrightarrow \infty} E Y_{n}^{+} =: l_{1} > 0$. Then $\Lambda_{n}'(a)=(E Y_{n} e^{a Y_{n}})/ E e^{a Y_{n}}$ satisfies:
\begin{equation}
\label{Defdomain2}
\liminf_{n\longrightarrow \infty} \Lambda_{n}'(a) \geq l_{1}(e^{a l_{1}}-1)/ M_{1} > 0,
\end{equation}
\noindent
where $M_{1} < +\infty $ is a bound of $E e^{a Y_{n}}$ (lemma \ref{lmLB2}), so that (\ref{Defdomain1}) is proved.
\\

By definition, (\ref{Defdomain1}) implies that for any $\varepsilon > 0$, there is a $n_{1}$ such that for all $n\geq n_{1}$, we have $ \Lambda_{n}'(a) > l-\varepsilon$; 
by taking $\varepsilon = l/2$ we have  $ \Lambda_{n}'(a) \geq l/2$ for all $n \geq n_{1}$.
\\
This means that the  interval $[0,l/2] \subset {\rm Im} \Lambda_{n}' (=\Lambda_{n}'(\R))  $ if $n \geq n_{1}$, and therefore $\Lambda^{*'}_{n} =\Lambda^{' -1}_{n}$ is defined on  $[0,l/2]$ for all $n \geq n_{1}$.
\\
Now for all $x \in [0,l/2]$ we have  ${\Lambda^{*}_{n}}'(x)= {\Lambda'}^{-1}_{n}(x) \leq a \leq K $ which implies that:
\begin{equation}
\label{BoundLambda*}
\Lambda^{*}_{n}(x) \leq K x \; \; \forall x \in [0,l/2], \forall n\geq n_{1}.
\end{equation}
\noindent
Hence, $\Lambda^{*}(x):=\limsup_{n\longrightarrow \infty} \Lambda^{*}_{n}(x) \leq K x, \forall x\in[0,b] $ with $b=l/2$, we note that $b$ may depends on $K$. This completes the proof of Proposition \ref{PropLB6} $\Box$
\\
\\ 
\noindent
{\bf Proof of Proposition \ref{PropLB5}}\\
Let $h \in \R$ and $C_{h} = \{ x \in \R: x \geq h \}$. We start from the inequality (\ref{LBineq}) and we consider the probability measures $Q_{N,h}$ such that $E(Q_{N,h}) \in C_{h}$, where for a probability measure $P$, $E(P)$ denotes the expectation or resultant of $P$, see \cite{Csiszar}.
\noindent
Then we have:
\begin{eqnarray*}
\frac{1}{N} \log P_{N}^{(N)} ( A_{N}) &\geq & - D( Q_{N,h} || P_{N}) +  \frac{1}{N}T_{N}{(1)}' + \frac{1}{N}T_{N}{(3)}'  \\
                                      &\geq & -\inf_{E(Q_{N,h}) \in C_{h}} D( Q_{N,h} || P_{N}) +  \frac{1}{N}T_{N}{(1)}' + \frac{1}{N}T_{N}{(3)}'
\end{eqnarray*}
We shall focus on the first term and the aim is to get an expression of 
$$ \inf_{E(Q_{N,h}) \in C_{h}} D( Q_{N,h} || P_{N})$$ 
\noindent
and its limit as $N\longrightarrow \infty$. From now on we specialize the proof to the case related to our initial problem: $-X_{N,i}$ are real random variables with the same law $P_{N}$, and $Q_{N,h}$ will be the law of $-X_{N,i} + h $. Then it is easy to see that $P_{N}\equiv Q_{N,h} $. By Csisz\'ar \cite{Csiszar}, Theorem 2 (or Theorem 3 which gives the same result here) we have:
\begin{equation}
D( C_{h} || P_{N}) = D( Q_{N,h}^{*} || P_{N}) = \inf_{E(Q_{N,h}) \in C_{h}} D( Q_{N,h} || P_{N}) = \max_{l\in \R^{+}} [ lh - \log E^{P_{N}} e^{-l X_{N,i}} ]
\end{equation}
This is nothing but the Legendre-Fenchel transform $\Lambda_{n}^{*}(h)$ of the $P_{N}$ or $-X_{N,i}$. With the notation $\Lambda_{n}(h)= \log E^{P_{N}} e^{-h X_{N,i}} $ we have indeed:
$$ \Lambda_{n}^{*}(h) = \max_{l\in \R} [ lh - \Lambda_{n}(l)] = \max_{l\in \R^{+}} [ lh - \Lambda_{n}(l)]  $$
\noindent
because $EX_{N,i}=0$. Then
\begin{equation}
\label{LB7}
\frac{1}{N} \log P_{N}^{(N)} ( A_{N}) \geq  - \Lambda_{n}^{*}(h) +  \frac{1}{N}T_{N}{(1)}' + \frac{1}{N}T_{N}{(3)}'
\end{equation}
and:
\begin{equation}
\label{LB8}
\liminf_{N\longrightarrow \infty}\frac{1}{N} \log P_{N}^{(N)} ( A_{N}) \geq  -\limsup_{N\longrightarrow \infty} \Lambda_{n}^{*}(h) + \liminf_{N\longrightarrow \infty} \frac{1}{N}T_{N}{(1)}' + \liminf_{N\longrightarrow \infty} \frac{1}{N}T_{N}{(3)}'
\end{equation}
\noindent
The liminf of the last two terms of (\ref{LB7}) being $\geq 0$ by the proof of Proposition \ref{PropLB2}, and with $\Lambda^{*}(h):=\limsup_{N\longrightarrow \infty} \Lambda_{n}^{*}(h)$ we get finally:
\begin{equation}
\label{LB9}
\liminf_{N\longrightarrow \infty}\frac{1}{N} \log P_{N}^{(N)} ( A_{N}) \geq  -\Lambda^{*}(h)
\end{equation}
\noindent
According to Proposition \ref{PropLB6}, $\Lambda^{*}$ is well defined on a set $[0,b]$ where it is finite and moreover for a given constant $K>0$, we can choose $b$ such that $\Lambda^{*}(h) \leq K h, \forall h \in [0,b] $. This completes the proof of Proposition \ref{PropLB5}. $\Box$
\subsection{The Laplace Method and large deviations}

After the transformation $\psi_{n}(x)= \varphi_{n}(x) / \sqrt{c_{n}}$, the action may be written as:
\begin{align*}
\label{eq441}
{\cal A}_{n} (\varphi) & = \int_{V}l_{n}^{(1})(\varphi_{n})(x)dx \\
             & =  u_{n}\int_{V}l_{n}(\varphi_{n})(x)dx\\
						 & =  u_{n} {\cal A}_{n}^{(1)} (\varphi)
\end{align*}

\noindent
where, in one of the cases that will be distinguished below, $u_{n}=g_{n}c_{n}^{2} \longrightarrow\infty$ and ${\cal A}_{n}^{(1)}(\varphi) \longrightarrow 0 $ a.e., obeys to a kind of law of large numbers. As previously mentioned ( \S 4.1), the investigation regarding the limiting field leads to the study of the estimate of:

\begin{equation}
\label{Laplace1}
Z_{n,V} = E e^{-{\cal L}_{n}^{R}} = E e^{u_{n} (- {\cal A}_{n}^{(1)} (\varphi))} = \int_{B} e^{ N F(\varphi)} \mu_{N} (d\varphi)
\end{equation}

\noindent
with $N=u_{n}$. Such an estimate are performed with generalizations of the Laplace method and its various generalizations in infinite dimension ; we refer to Pitebarg and Fatalov \cite{PF} for a detailed review of this topic. In large deviations contexts, the Varadhan lemma is often used and may be formulated in fairly general settings. Roughly speaking, the Laplace method tells us that if a sequence of measures $\nu_{N}$ on some space $B$ satisfies:
\begin{equation}
\label{Laplace2}
d \nu_{N} (\varphi) \sim  e^{-N I (\varphi)} \; {"} d\varphi {"} 
\end{equation}
with $I:B\rightarrow R$ being some {\it rate function}, then, for a function $F:B\rightarrow R$ having some properties, we will have:
\begin{equation}
\label{Laplace3}
Z_{N} = E_{\nu_{N}} e^{N F} = \int_{B} e^{ N F(\varphi)} \nu_{N} (d\varphi) \sim e^{N \sup_{\varphi}(F(\varphi) - F(\varphi)) }
\end{equation}
Now, the condition (\ref{Laplace2}) is correctly formulated by requiring that the sequence of measures $\nu_{N}$ satisfies a large deviation principle (LDP) with a rate function $I$: 
\begin{equation}
\label{Laplace4}
-\inf_{{\rm int} A} I(x)    \leq \lim\inf_{N\rightarrow +\infty}    \frac{1}{N} \log \nu_{N} (A) 
               \leq \lim\sup_{N\rightarrow +\infty}  \frac{1}{N} \log \nu_{N} (A) \leq -\inf_{{\rm cl} A} I(x), 
\end{equation}
for all $A\in {\cal B}$ (the Borel $\sigma$-field of $B$) ; ${\rm int}$ and ${\rm cl} A$ are the interior and closure of the set $A$ w.r.t the topology of $B$ .
The rate function $I: B \longrightarrow [0, +\infty]$ is assumed to be lower semi-continuous. If in addition the
level sets $\{ x\in E: I(x)\leq L\}, L\geq 0$ are compact, then $I$ is  said to be a good or proper rate function. The sequence $\nu_{N}$ may be a family $P_{\varepsilon}$  depending on $\varepsilon \rightarrow 0 $ ($\varepsilon=1/N)$, in which case (\ref{Laplace4}) is written as
\begin{equation}
\label{Laplace5}
-\inf_{{\rm int} A} I(x)\leq \lim\inf{\varepsilon\rightarrow 0}
 \varepsilon\log P^{\varepsilon}(A) \leq \lim\sup_{\varepsilon\rightarrow 0}
  \varepsilon\log P^{\varepsilon}(A)\leq -\inf_{{\rm cl} A} I(x).
\end{equation}
Although the measure space $B$ is often supposed to be metric, separable and complete (Polish space), many large deviations results are valid in more general settings (and this may be actually needed, cf., e.g., \cite{AA1} where the reference space $D([0,1], \R)$, the set of piecewise continuous functions is non separable when endowed with the supremum norm $\|.\|_{\infty}$). In our case, the reference probability space is ${\cal P}={\cal S}^{'}(\R^{d})$, but we shall work with the arrays of real random variables $X_{ni}$.  
\\
\noindent
Now, when the sequence of measures $\nu_{N}$ satisfies the LDP (\ref{Laplace4}), the estimate (\ref{Laplace3}) is formulated as follows: for every bounded function $F \longrightarrow\R$, we have:
\\
\begin{equation}
\label{Laplace2}
 \int_{E} e^{n F(x)} \mu_{n}(dx) = e^{n (\sup_{x \in E}(F(x)-I(x))+o(1))},
\end{equation}
or:
\begin{equation}
\label{Laplace3}
 \lim_{n \longrightarrow\infty} \frac{1}{n}\log \int_{E} e^{n F(x)} \mu_{n}(dx)= \sup_{x \in E}(F(x)-I(x)).
\end{equation}
\\
The Varadhan lemma is usually stated with assumptions that are not be fulfilled in our case. As we need only the lower bound part of this lemma, we recall it in a following form:
\begin{thm}\label{VaradhanLm}
Let $\mu_{n}$ be a sequence of probability measures on a Polish space $E$ and $F: E \longrightarrow \R$ a function.
\\
(1) Lower bound: Suppose that $\mu_{n}$ has a large deviations lower bound, that is: there is a function $I: E \longrightarrow \R$, such that for each open set $C \subset E$:
\begin{equation}
\liminf_{n \rightarrow \infty} \frac{1}{n}\log \mu_{n}(C) \geq -\inf_{x\in A} I(x) 
\end{equation}
Then, if $F$ is lower semi-continuous, we have:
\begin{equation}
\label{VaradhanLB}
 \lim_{n \longrightarrow\infty} \frac{1}{n}\log  \int_{E} e^{ n F(x)} \mu_{n}(dx) \geq \sup_{x \in E}(F(x)-I(x)).
\end{equation}
\noindent
(2) Upper bound: Suppose that $\mu_{n}$ satisfies a large deviations upper bound, that is, there is a proper rate function $I: E \longrightarrow \R$ (which is lower semi-continuous and its level sets $\{x: I(x) \leq r\}, r\in\R$ are compact) such that for each closed set $C \subset E$:
\begin{equation}
\limsup_{n \rightarrow \infty} \frac{1}{n} \log \mu_{n}(C) \leq -\inf_{x \in A} I(x) 
\end{equation}
Then, if $F$ is lower semi-continuous and bounded above, we have:
\begin{equation}
\label{VaradhanUB}
 \limsup_{n \longrightarrow\infty} \frac{1}{n}\log \int_{E} e^{ n F(x)} \mu_{n}(dx) \leq  \sup_{x \in E}(F(x)-I(x)).
\end{equation}
And therefore if the conditions of (1) and (2) are fulfilled, we have:
\begin{equation}
\label{VaradhanLemma}
 \lim_{n \longrightarrow\infty} \frac{1}{n}\log  \int_{E} e^{ n F(x)} \mu_{n}(dx) = \sup_{x \in E}(F(x)-I(x)).
\end{equation}
\end{thm}
Let us remark that for the lower bound part, no assumption is made about $I$, and the only condition set for $F$ is the
lower semi-continuity. Its proof is quite simple and it is the only part used in this paper ; we remind it here:
\\
By the lower semi-continuity of $F$, for each $\varepsilon > 0$ and $x_{0} \in E$, there exists a neighborhood $B_{x_{0}}$ of $x_{0}$ such that
$F(x) \geq F(x_{0}) - \varepsilon$ for all $x \in B_{x_{0}}$ ; so that we have:
\begin{align*}
\nonumber
Z_{n}= \int_{E} e^{n F(x)} \mu_{n}(dx)&\geq \int_{B_{x_{0}}} e^{n F(x)} \mu_{n}(dx) \\
                                        &\geq  \int_{B_{x_{0}}} e^{n (F(x_{0})-\varepsilon)} \mu_{n}(dx) \\
                                         &= e^{ n (F(x_{0})-\varepsilon)} \mu_{n}(B_{x_{0}}) 
\end{align*}
and
\[ \frac{1}{n}\log Z_{n} \geq F(x_{0})-\varepsilon + \frac{1}{n}\log \mu_{n}(B_{x_{0}})   \]

\noindent
By the lower bound assumption on the $\mu_{n}$ we get:

\[ \liminf_{n \rightarrow \infty} \frac{1}{n}\log Z_{n} \geq F(x_{0})-\varepsilon -\inf_{x\in B_{x_{0}} } I(x) \geq F(x_{0})-\varepsilon - I(x_{0})   \]
This inequality being true for all $\varepsilon > 0$ and $x_{0}$, we get the lower bound (\ref{VaradhanLB}). $\Box$


\subsection{End of the proof of the main theorem}

The action is written in the following form:
\begin{align}
\label{proof1}
 {\cal A}_{n} &=  g_{n} c_{n}^{2} \int_{V} :\psi_{n}^{4}(x): dx + m_{n}c_{n} \int_{V} :\psi_{n}^{2}(x): dx + a_{n} d_{n }c_{n} \frac{1}{d_{n}}\int_{V}:(\partial \psi_{n})^{2}(x): dx \\ \nonumber
                  &= g_{n}'\int_{V} :\psi_{n}^{4}(x): dx + m_{n}'\int_{V} :\psi_{n}^{2}(x): dx + a_{n}' \frac{1}{d_{n}}\int_{V}:(\partial \psi_{n})^{2}(x): dx
\end{align}

\noindent
where $d_{n }$ is chosen in order to ensure that ${\rm E}((1/d_{n})\int_{V}:(\partial \psi_{n})^{2}(x): dx)^{2}$ converges to some $K < \infty$. We distinguish the two cases of the theorem:
\\
\\
\noindent
$\bullet$ {\bf Case (A):} At least one of the terms $g_{n} c_{n}^{2}$ or $m_{n}c_{n}$ or $a_{n} d_{n }c_{n}$ is not bounded: 
\\
\\
\noindent
Following the remark of $\S 3$ we may suppose the limit of the unbounded term(s) is $+\infty$. We shall discuss 3 cases, depending of the dominating sequence of the above-mentioned 3 terms $g_{n}', m_{n}'$ or $a_{n}$:
\\
\\
\noindent
$\bullet$ {\bf Case (A.1):}  $g_{n}'=g_{n} c_{n}^{2}$ is the dominating factor in (\ref{proof1}) and $g_{n}' \longrightarrow \infty $ and $\alpha_{n}^{0} = O(1)$ and $\beta_{n}^{0}= O(1)$ with 
\[ \alpha_{n}^{0} = \frac{m_{n}}{c_{n} g_{n}} ,  \; \; \beta_{n}^{0} = \frac{a_{n} d_{n }}{c_{n} g_{n}} \]
The previous two sequences are bounded and as we can consider subsequences, we may suppose that they converge respectively to some $\alpha$ and $\beta$. We also might have $\alpha=0$ or $\beta=0$ or both.

\noindent
Let us rewrite the modified action and the exponent ${\cal A}_{n}^{(1)}$ as:
\begin{eqnarray*}
 {\cal A}_{n} &=&  g_{n} c_{n}^{2} [\int_{V} :\psi_{n}^{4}(x): dx - \frac{m_{n}}{c_{n} g_{n}} \int_{V} :\psi_{n}^{2}(x): dx - \frac{a_{n} d_{n }}{c_{n} g_{n}} \frac{1}{d_{n}}\int_{V}:(\partial \psi_{n})^{2}(x): dx \\
                  &=:& g_{n} c_{n}^{2} {\cal A}_{n}^{(1)}(\varphi),
\end{eqnarray*}
And ${\cal A}_{n}^{(1)}(\varphi)$ is of the form:
\begin{equation}
\label{LagCase1}
{\cal A}_{n}^{(1)}(\varphi)=  \lambda_{n} \int_{V} :\psi_{n}^{4}(x): dx - \alpha_{n}\int_{V} :\psi_{n}^{2}(x): dx - \beta_{n} \frac{1}{d_{n}}\int_{V}:(\partial \psi_{n})^{2}(x): dx,
\end{equation}
\noindent
where with $\lambda_{n} \equiv 1, \alpha_{n} \longrightarrow \alpha$ and $\beta_{n} \longrightarrow  \beta $. We can apply the results obtained in the previous sections: 
${\cal A}_{n}^{(1)}(\varphi) $ may be transformed as a mean of an array of random variables, which makes the link with the possibility of using large deviations techniques and other probabilistic results. By the law of large numbers (Proposition \ref{PropLLN}):
\begin{equation}
\label{Proof10}
{\cal A}_{n}^{(1)}(\varphi) = 1/n^{d}\sum_{i=1}^{n^{d}} X_{n,i} \longrightarrow 0, \; a.e. 
\end{equation}
Let:
\begin{equation}
R_{n}= \frac{d\mu_{n,V}}{d\mu_{0}}= \frac{1}{Z_{n,V}}e^{-g_{n} c_{n}^{2} {\cal A}_{n}^{(1)}(\varphi)} \; \;  Z_{n,V}=  {\rm E} e^{-g_{n} c_{n}^{2} {\cal A}_{n}^{(1)}(\varphi)}
\end{equation}
At this point we make the assumption that $g_{n} n^{d-4}$ does not converge to $0$. This means that in dimension 4 we suppose that $g_{n}$ does not converge to $0$, and in dimensions $d\geq 5$ we may accept that $g_{n}$ converge to $0$ but with a speed less that $1/n^{d-4}$: 
\\
\\
{\bf (A.1.1):} $g_{n} n^{d-4}$ does not converge to $0$, i.e. $\limsup_{n} g_{n} n^{d-4} = G_{1} > 0 $
\\
\\
This assumption is legitimate only for dimension $d\geq 4$: in lower dimensions we can not have  $g_{n} n^{d-4} \longrightarrow 0$ when the coupling $g_{n}$ is constant for $d\leq 3$. We shall also see that the following arguments do not work for a negative coupling constant sequence.
\noindent
As we can consider subsequences, we may suppose with this assumption that there is $N_{1} \geq 0$ such that for all $n\geq N_{1}$ : $  g_{n} n^{d-4} \geq G > 0 $
\\
On the other hand $c^{2}_{n} \sim K n^{2(d-2)}$, and for large $n$ ($n\geq N_{2}$ for some $N_{2}$)  we have: 
\[    (K/2) n^{2(d-2)} \leq c^{2}_{n} \leq 2K n^{2(d-2)}  \] 
and then:
\begin{equation}
\label{Proof12}
g_{n} c^{2}_{n} \geq \frac{G K }{2} n^{d}  \; \; {\rm for} \; n\geq N_{3}:= \max({N_{1},N_{2}})
\end{equation}
Therefore:
\begin{align}
\label{Proof}
{\rm E} e^{- g_{n} c_{n}^{2} {\cal A}_{n}^{(1)}} & \geq {\rm E} \exp (- g_{n} c_{n}^{2} {\cal A}_{n}^{(1)})1_{{\cal A}_{n}^{(1)}\leq 0} \nonumber \\
                             & \geq  {\rm E} \exp ( - \frac{G K }{2} n^{d} {\cal A}_{n}^{(1)}) 1_{{\cal A}_{n}^{(1)}\leq 0} \nonumber \\ 
                             & \geq  {\rm E} \exp ( - \frac{G K }{2} n^{d} {\cal A}_{n}^{(1)}) -1
\end{align}
and 
\begin{equation}
\label{Proof14}
\frac{1}{n^{d}} \log Z_{n,V} = \frac{1}{n} \log {\rm E} e^{- g_{n} c_{n}^{2} {\cal A}_{n}^{(1)}} \geq \frac{1}{n^{d}} \log {\rm E} \exp ( - \frac{G K }{2} n^{d} {\cal A}_{n}^{(1)}) -\frac{1}{n^{d}}
\end{equation}
(we use $\log(x-1) \geq \log x -1$ for $x\geq 2$, the corresponding $x$ in the last equation is indeed large (hence $\geq 2$) as we shall see.). This implies that:

\begin{align}
\label{Proof15}
\liminf_{n\rightarrow +\infty}\frac{1}{n^{d}} \log Z_{n,V} & 
                                          \geq \liminf_{n\rightarrow +\infty} \frac{1}{n^{d}} \log {\rm E} e^{- g_{n} c_{n}^{2} {\cal A}_{n}^{(1)}} \\ \nonumber
                                     & \geq \liminf_{n\rightarrow +\infty} \frac{1}{n^{d}} \log {\rm E} \exp ( - \frac{G K }{2} n^{d} {\cal A}_{n}^{(1)}) -0
\end{align}

(1) By the lower bound result (Proposition \ref{PropLB5}) for the array $X_{n,i}$ or the sequence ${\cal A}_{n}^{(1)}$ given by (\ref{LagCase1}), and with $N=n^{d}$, there exists a function $\Lambda^{*}$ such that, for any open set $C\subset \R$:
\begin{equation}
\label{Proof16}
\liminf_{N\longrightarrow \infty}\frac{1}{N} \log {\rm Pr} ( {\cal A}_{n}^{(1)} = 1/n^{d}\sum_{n=1}^{n^{d}} X_{n,i} \in C ) \geq  -\inf_{h\in C}\Lambda^{*}(h)
\end{equation}
\noindent
with the properties that $\Lambda^{*}$ is well defined and is finite on an interval $[0,b], b >0$ and furthermore $b$ can be chosen so that $\Lambda^{*}(h) \leq G K h/4, \forall h\in [0,b]$.
\\
\\
(2) We can the use the lower bound of the Varadhan lemma, applied in fact to the sequence $-{\cal A}_{n}^{(1)}$:
\[ \liminf_{N\longrightarrow \infty} \frac{1}{N} \log {\rm E} e^{N F(-{\cal A}_{n}^{(1)})} \geq \sup_{h\in \R} [F(h) - \Lambda^{*}(h) ] \]

\noindent
We apply this formula with $F(h)=G K h/2, N=n^{d}$ (notice that $F$ is not bounded) and we get:
\begin{align}
\label{Proof17}
\liminf_{n\longrightarrow \infty} \frac{1}{n^{d}} \log {\rm E} \exp(-\frac{G K}{2}n^{d}{\cal A}_{n}^{(1)}) &\geq  \sup_{h\in [0,b} [G K h/2  - \Lambda^{*}(h)] \\ \nonumber
                                                                        &\geq \sup_{h\in [0,b]} [K h/2 - \frac{G K h}{4}] \\  \nonumber
																																         &\geq \frac{G K b}{4},
\end{align}
This, with (\ref{Proof15}) shows that: $\liminf_{n} Z_{n,V} = +\infty$ and we have even more:
\begin{equation}
\label{Proof18}
\liminf_{n\rightarrow +\infty}\frac{1}{n^{d}} \log Z_{n,V} \geq z_{\infty}:= \frac{G K b}{4} > 0
\end{equation}

(3)We turn now to the limit of the Radon-Nikodym density of the interacting field. Since $  g_{n} n^{d-4} \geq G > 0 $ and $(K/2) n^{2(d-2)} \leq c^{2}_{n} \leq 2K n^{2(d-2)} $ for $N\geq N_{2}$ , we write this time:
\[   g_{n} c_{n}^{2} \geq \frac{K}{2}  g_{n} n^{2d-4} = \frac{K G n^{d} }{2} \frac{g_{n}n^{d-4} }{G}   \]
and since $g_{n}n^{d-4}/G \geq 1 $ we have:
\[ \E \exp (-g_{n} c_{n}^{2} {\cal A}_{n}^{(1)}) \geq ( \E \exp(-\frac{K G n^{d} }{2}{\cal A}_{n}^{(1)}) )^{g_{n}n^{d-4}/G}, \]
which implies that:
\begin{align}
\label{Proof19}
\log R_{n} & \leq g_{n} c_{n}^{2} [ -{\cal A}_{n}^{(1)} -  \frac{g_{n} n^{d-4}}{G g_{n}c_{n}^{2}} \log {\rm E} \exp(-\frac{K G n^{d} }{2}{\cal A}_{n}^{(1)})] \\ \nonumber
           & \leq g_{n} c_{n}^{2} [ -{\cal A}_{n}^{(1)} -  \frac{g_{n} n^{d-4}}{G g_{n} 2 K n^{2d-4} } (\log {\rm E} \exp(-\frac{K G n^{d} }{2}{\cal A}_{n}^{(1)}) ] 
\end{align}
where we have used $c^{2}_{n} \leq 2K n^{2(d-2)} $ for $N\geq N_{2}$. Hence:
\begin{align}
\label{Proof110}
\limsup_{n\longrightarrow \infty} \frac{1}{g_{n} c_{n}^{2}}  \log R_{n} & \leq 
                      \limsup_{n\longrightarrow \infty} [ -{\cal A}_{n}^{(1)} -  \frac{1}{2 K G n^{d}} \log {\rm E} \exp(-\frac{K G n^{d} }{2}{\cal A}_{n}^{(1)}) ]\\ \nonumber
           & \leq \limsup_{n\longrightarrow \infty} (-{\cal A}_{n}^{(1)}) -\liminf_{n\longrightarrow \infty} \frac{1}{2 K G n^{d}} \log {\rm E} \exp(-\frac{K G n^{d} }{2}{\cal A}_{n}^{(1)})
\end{align}
By the law of large numbers (\ref{Proof10}) and the lower bound (\ref{Proof17}) we get:
\begin{equation}
\label{Proof111}
\limsup_{n\longrightarrow \infty} \frac{1}{g_{n} c_{n}^{2}}  \log R_{n} \leq -\frac{b}{4}
\end{equation}
which implies that: 
\begin{equation}
\label{Proof112}
\limsup_{n\longrightarrow \infty} \frac{d\mu_{n,V}}{d\mu_{0}}= 0, \; {\rm a.e.}
\end{equation}
\noindent
This completes the proof of Theorem 3.1 (Part (2)) in the case {\bf( A.1)}.
\\
\\
\noindent
$\bullet$ {\bf Case (A.2): } The dominating term is  $m_{n} c_{n} \longrightarrow \infty$: That is 

\[ \frac{g_{n}{c_{n}}}{ m_{n}} = O(1)  \; {\rm and } \; \frac{a_{n} d_{n }}{m_{n}}= O(1)   \]
In this case the action and the exponent ${\cal A}_{n}^{(1)}$ can be written as:
\begin{eqnarray*}
 {\cal L}^{R}_{n} &=&  m_{n} c_{n} [ \frac{g_{n}{c_{n}}}{ m_{n}}  \int_{V} :\psi_{n}^{4}(x): dx - \int_{V} :\psi_{n}^{2}(x): dx - 
                       \frac{a_{n} d_{n }}{m_{n}} \frac{1}{d_{n}}\int_{V}:(\partial \psi_{n})^{2}(x): dx] \\
                  &=& m_{n} c_{n} {\cal A}_{n}^{(1)}
\end{eqnarray*}
So that  ${\cal A}_{n}^{(1)}$ is of the form:
\begin{equation}
{\cal A}_{n}^{(1)}(\varphi)=  \lambda_{n} \int_{V} :\psi_{n}^{4}(x): dx - \alpha_{n}\int_{V} :\psi_{n}^{2}(x): dx - \beta_{n} \frac{1}{d_{n}}\int_{V}:(\partial \psi_{n})^{2}(x): dx,
\end{equation}
\noindent
where we may suppose that $\lambda_{n}$ and $\beta_{n}$ have a limit (possibly $=0$) and $\alpha_{n}\equiv 1$ ; ${\cal A}_{n}^{(1)}$ may be written as the mean of an array of random variables. We use as before Proposition \ref{PropLB5} to the get a large deviations lower bound for the $X_{n,i}$; we conclude the proof exactly as in the case A.1 provided we add the assumption:
\\
\\
{\bf Condition (A.2.1):} $m_{n} c_{n} \geq K_{4} n^{d}$ for $n$ sufficiently large and for some constant $K_{4}$. Or at least this inequality holds for a subsequence of $n$.
\\
\\
We note that:
\\
$\bullet$ Condition (A.2.1) is automatically satisfied if (A.1.1) is satisfied (i.e. in dimension 4 that $g_{n}$ does not tend to $0$).
\\
Indeed In this case (A.2), we have for some constant $K'$: $m_{n} c_{n} \geq K' g_{n}c_{n}^{2} \geq K' g_{n} (K/2) n^{2d-4} \geq K' (K/2) G n^{d} $. 
\\
In the case where (A.1.1) is not satisfied i.e. $g_{n}n^{d-4} \longrightarrow 0$, condition (A.2.1) means that $ m_{n} \geq K_{4} n^{2}$ for some constant $K_{4}$.
\\
\\
\noindent
$\bullet${\bf Case (A.3): } The dominating term is  $a_{n} c_{n} d_{n} \longrightarrow \infty$: That is 

\[ \frac{g_{n}{c_{n}}}{a_{n} d_{n}} = O(1)  \; {\rm and } \; \frac{m_{n}}{ a_{n}d_{n}}  = O(1)   \]
\noindent
In this case we write the action and the exponent ${\cal A}_{n}^{(1)}$ are written as:
\begin{align*}
 {\cal L}^{R}_{n} &=  a_{n} c_{n} d_{n} [ \frac{g_{n}c_{n}}{ a_{n}d_{n}}  \int_{V} :\psi_{n}^{4}(x): dx - \frac{m_{n}}{ a_{n}d_{n}} \int_{V} :\psi_{n}^{2}(x): dx - 
                       \frac{1}{d_{n}}\int_{V}:(\partial \psi_{n})^{2}(x): dx] \\
                  &=  a_{n} c_{n} d_{n} {\cal A}_{n}^{(1)}
\end{align*}
So that  ${\cal A}_{n}^{(1)}$ is of the form:
\begin{equation}
{\cal A}_{n}^{(1)}(\psi)=  \lambda_{n} \int_{V} :\psi_{n}^{4}(x): dx - \alpha_{n}\int_{V} :\psi_{n}^{2}(x): dx - \beta_{n} \frac{1}{d_{n}}\int_{V}:(\partial \psi_{n})^{2}(x): dx,
\end{equation}
\noindent
where we may suppose this time that $\lambda_{n}$ and $\alpha_{n}$ have a limit (possibly $=0$) and $\beta_{n}=1$. ${\cal A}_{n}^{(1)}$ may be written as the mean of an array of random variables in the same way as the case $N^{\circ} A.1 $. The assertions of Proposition \ref{PropLB5} related to the lower bound can be applied. And with the same arguments as in the case A.1 we will get $\limsup_{n} R_{n}=0 \;$ a.e. and $\limsup_{n} Z_{n,V}=+\infty$, provided we add the assumption:
\\
\\
{\bf Condition (A.3.1):} $a_{n} c_{n} d_{n} \geq K_{5} n^{d}$ for $n$ sufficiently large and for some constant $K_{5}$. Or at least this inequality holds for a subsequence of $n$.
\\
\\
We note that:  
\\
$\bullet$ Condition (A.3.1) is automatically satisfied if (A.1.1) is satisfied (i.e. in dimension 4 that $g_{n}$ does not tend to $0$).
\\
Indeed, in this case (A.3) we have for some constant $K'$:
\[  a_{n} c_{n} d_{n} \geq K' g_{n}c_{n}^{2} \geq K' g_{n} (K/2) n^{2d-4} \geq K' (K/2) G n^{d}. \]
We recall $d_{n}\sim K" n^{2}$. Then in the case where (A.1.1) is not satisfied i.e. $g_{n}n^{d-4} \longrightarrow 0$, condition (A.3.1) means that $ a_{n}n^{d-2}n^{2} \geq K_{6} n^{d}$ for some constant $K_{6}$, i.e. $\limsup_{n} a_{n} > 0$.
\\
\\
\noindent
{\it Remark.} The above arguments can be summarized as follows: The Radon-Nikodym derivative $R_{n}$ can be written as
\begin{equation}
R_{n}=\frac{d\mu_{n}}{d\mu_{0}}= \frac{1}{Z_{n}}{ \exp (-g_{n} c_{n}^{2}{\cal A}_{n}^{(1)})}
            = \exp [g_{n} c_{n}^{2}(-{\cal A}_{n}^{(1)}) - \frac{1}{g_{n} c_{n}^{2}}\log \E e^{-g_{n} c_{n}^{2}{\cal A}_{n}^{(1)}} ]
\end{equation}
If $g_{n}n^{d-4}\nrightarrow 0$, i.e. $g_{n}n^{d-4} \geq K_{1}> 0$ for some $K_{1}$ (we can suppose this as being valid for every $n$), then:
\[ \frac{g_{n} c_{n}^{2} }{n^{d}} \sim \frac{K' g_{n} n^{2(d-2) }}{n^{d}}  \sim K' g_{n} n^{d-4} \geq K' K_{1} = K_{2} > 0\]
and:
\begin{align*}
\E \exp [-g_{n} c_{n}^{2}{\cal A}_{n}^{(1)}]& = \E \exp [ \frac{g_{n} c_{n}^{2}}{K_{2}n^{d}} K_{2}n^{d}(-{\cal A}_{n}^{(1)})] \\
               & \geq (\E \exp [ K_{2}n^{d}(-{\cal A}_{n}^{(1)})] )^{\frac{g_{n} c_{n}^{2}}{K_{2}n^{d}}}
\end{align*}
because $g_{n} c_{n}^{2} / (K_{2}n^{d}) \geq 1$. This yields:
\begin{equation*}
\frac{1}{g_{n} c_{n}^{2}}\log \E e^{-g_{n} c_{n}^{2}{\cal A}_{n}^{(1)}} \geq \frac{1}{K_{2}n^{d}} \log \E e^{K_{2}n^{d}(-{\cal A}_{n}^{(1)})}.
\end{equation*}
and by the lower bound of large deviations of ${\cal A}_{n}^{(1)}$ and Varadhan's lemma we get:
\begin{equation}
\label{Proof200}
\liminf{n\rightarrow \infty}\frac{1}{g_{n} c_{n}^{2}}\log \E e^{-g_{n} c_{n}^{2}{\cal A}_{n}^{(1)}} \geq \frac{1}{K_{2}} \max[K_{2} h - \Lambda^{*}(h) ].
\end{equation}
The properties of $\Lambda^{*}(h)$ mentioned in propositions \ref{PropLB5} and \ref{PropLB5} imply that $\Lambda^{*}(h)\leq K_{2} h/2$ for $h\in [0,b]$ with $b>0$ sufficiently small and we deduce from (\ref{Proof200}) that:
\begin{equation*}
\liminf_{n\rightarrow \infty}\frac{1}{g_{n} c_{n}^{2}}\log \E e^{-g_{n} c_{n}^{2}{\cal A}_{n}^{(1)}} \geq \frac{b}{2}
\end{equation*}
and
\begin{equation}
\liminf_{n\rightarrow \infty}\frac{1}{g_{n} c_{n}^{2}}\log R_{n} \leq 0 - \liminf{n\rightarrow \infty}\frac{1}{g_{n} c_{n}^{2}}\log \E e^{-g_{n} c_{n}^{2}{\cal A}_{n}^{(1)}}  \leq - \frac{b}{2}
\end{equation}
which means that $ R_{n} \rightarrow 0$, a.e. 
\\
\\
\noindent
$\bullet$ {\bf Case B: } The terms $g_{n} c_{n}^{2}, m_{n}c_{n}$ and $a_{n} d_{n }c_{n}$ are bounded: \\
\\
In this case, since we have:
\[  \lim_{n\longrightarrow \infty} \int_{V} :\psi_{n}^{4}(x): dx = \lim_{n\longrightarrow \infty}  \int_{V} :\psi_{n}^{2}(x): dx = \lim_{n\longrightarrow \infty}  
                       \frac{1}{d_{n}}\int_{V}:(\partial \psi_{n})^{2}(x): dx = 0 \]
almost everywhere, we get $\lim_{n\longrightarrow \infty}  {\cal A}^{(1)}_{n} = 0, \; {\rm a.e.}  $, and therefore:
\begin{equation}
R_{\infty} := \lim_{n\longrightarrow \infty} R_{n,V}= \lim_{n\longrightarrow \infty} \frac{d\mu_{n,V}}{d\mu_{0}}= 1, \; a.e.
\end{equation}
On the other hand, we have $\E R_{n,V} =1$ which obviously implies that $\lim_{n\longrightarrow \infty}  R_{n,V} = \E R_{\infty}$. 
\\
Now, it is a well known result that if a sequence of integrable random variables $X_{n}$ converges almost everywhere (or even in probability) to an integrable random variable $X$, and if $\lim_{n} \E |X_{n}| = \E |X|$, then $X_{n}$ converges to $X$ in $\L_{1}$ (cf., e.g., Neveu \cite{Neveu}, p.56, Ex. II-6-5). From this we deduce that $ R_{n,V}$ converges to $ R_{\infty}\equiv 1$ in $\L_{1} ({\cal S}'(\R^{d}))$, and therefore the sequence of measures $ \mu_{n,V}$ converges strongly (setwise) to $\mu_{0}$, which means the the limiting field $\varphi_{d}^{4}$ is the free field. This completes the proof of theorem 3.1.
\\
\noindent
Let us remark that the later arguments are valid for all dimensions $d\geq 1$ in the case where the terms $g_{n} c_{n}^{2}, m_{n}c_{n}$ and $a_{n} d_{n }c_{n}$ are bounded; while the former ones, corresponding to the case where at least one of the previous three terms is unbounded, are valid only for dimensions $d\geq 4$. 
\\
\\
\noindent
Finally, we also note that the scalar character of the field $\varphi$ (i.e. the fact that $\varphi(f) \in \R)$ does not play any role in
the intermediate propositions or the end of the proof and Theorem 3.1 is thus also valid for the $|\varphi|^{4}$ vector field.
$\Box$.

\newpage
{\small


\begin{thebibliography}{cc}
\bibitem{AA1} A. Aboulalaa, Grandes d\'eviations et principe de maximum d'entropie pour les processus de sauts, 
{\it C. R. Acad. Sc. Paris, S\'erie 1, Math\'ematiques}, 321 (7), 919-921, 1995.

\bibitem{AA2} A. Aboulalaa, The Gibbs principle of Markov jump processes, {\it Stoch. Proc. Appl.}, 64, 257-271, 1996

\bibitem{AA3} A. Aboulalaa, Random fields, large deviations and triviality in quantum field theory. Part II, arXiv preprint arXiv:2212.14879, 2022.

\bibitem{Aizenman} M. Aizenman, Geometric analysis of $\varphi_{4}^{4}$ fileds and Ising models, {\it Com. Math. Phys.}, 86 1-48, 1982, and
{\it Phys. Rev. lett.} 47, 1-4, 1981.

\bibitem{A-DC} M. Aizenman, and Hugo Duminil-Copin, Marginal triviality of the scaling limits of critical 4D Ising and $\varphi_{4}^{4}$ models, Annals of Mathematics 194.1 (2021): 163-235.

\bibitem{Araki} H. Araki, {\it Mathematical Theory of Quantum Fields }, Oxford Univ. Press, Oxford 1999. 

\bibitem{BR} R.R. Bahadur, M. Raghavachari, Some asymptotic properties of likelihood ratios in general sample spaces,
{Proc. Sixth Berkeley Symp. Probab.} 129-152, 1968.

\bibitem{BZG} R.R. Bahadur, S.L. Zabell, J.C. Gupta, Large deviations, tests and estimates, {\it in I.M. Chakravarti, ed. Asymptotic Theory
of Statistical Tests and Estimation, Proc. Adv. In. Symp., Chapel Hill, NC}, Academic Press Boston, pp. 33-64, 1980.

\bibitem{Balaban} T. Balaban, Ultraviolet stability in field theory. The $\varphi_{3}^{4}$ model,
In J. Fr\"ohlich (eds), {\it Scaling and self-Similarity in Physics}, Vol. 7, Birkh\"auser, Boston MA, (1983).

\bibitem{BF} G.A. Battle III, P. Federbush, A phase cell cluster expansion for a hierarchical $\varphi_{3}^{4}$ model, 
{\it Com. Math. Phys.}, 88 263-293, (1983).

\bibitem{Baumann} K. Baumann, When is a field theory a generalized free field?,  {\it Commun. Math. Phys.} 43, 221–223 (1975).

\bibitem{Benfatto-al} G. Benfatto, M. Cassandro, G. Gallavotti, F. Nicol\'o, E. Olivieri, E. Presutti, E. Scacciatelli,
Ultraviolet stability in Euclidean scalar field theories. {\it Commun. Math. Phys.} 71, 95-130, 1980

\bibitem{BFS} D.C. Brydges, J. Fr\"ohlich, A.D. Sokal, A new proof of the existence on nontriviality of the continum $\varphi_{2}^{4}$
and $\varphi_{3}^{4}$ quantum field theories, {\it Com. Math. Phys.}, 91 141-186, (1983)

\bibitem{BFS83} D.C. Brydges, J. Fr\"ohlich, A.D. Sokal, "The random-walk representation of classical spin systems and correlation inequalities." Communications in mathematical physics 91.1 (1983): 117-139.

\bibitem{Callaway} D. Callaway, Triviality pursuit: can elementary scalar particles exist ?, {\it Physics Reports} 
167, No 5, 241-320, 1988.

\bibitem{Collins} J.C. Collins, {\it Renormalization: an introduction to renormalization, the renormalization group and the operator-product expansion} Cambridge university press, 1985.

\bibitem{Csiszar} I. Csisz\'ar, Sanov property, generalized I-projection and a conditional limit theorem, 
{\it Ann. Probab.} 12 768-793, 1984;

\bibitem{Dimock} J. Dimock, {\it Quantum Mechanics and Quantum Field Theory}, Cambridge Univ. Press, 2013. 

\bibitem{FO} J.S. Feldman and K. Osterwalder, The Wightman axioms and the mass gap for weakly
coupled $\varphi_{3}^{4}$ quantum field theories, {\it Ann. Phys.}, 97, 80–135 (1976).

\bibitem{FFS} R. Fernandez, J. Fr\"ohlich, A.D.Sokal, Random Walks, {\it Critical Phenomena and Triviality in Quantum Field theory}, Springer Verlag, 1992.  

\bibitem{Fr74} J. Fr\"ohlich, Verification of axioms for Euclidean and relativistic fields and Haag’s theorem
in a class of $P(\varphi)_{2}$ models, {\it Ann. Inst. Henri Poincaré}, 21, 271–317 (1974).

\bibitem{Frohlich82} J. Fr\"ohlich, On the triviality of $\lambda \varphi_{3}^{4}$ theories and the approach to the critical point in $d\geq 4$
dimensions, {\it Nucl. Phys.}, B200, 281–296 (1982).

\bibitem{Gallavotti1} G. Gallavotti, Renormalization theory and ultraviolet stability for scalar fields via renormalization group methods, Reviews of Modern Physics, 57(2), 471 (1985).

\bibitem{Gallavotti2} G. Gallavotti, Constructive quantum field theory, arXiv preprint math-ph/0510014 (2005).


\bibitem{GR84} G. Gallavotti, V. Rivasseau,  $\varphi_{4}^{4}$ field theory in dimension 4: a modern introduction to its unsolved prolems,
 {\it Ann. Inst. Henri Poincaré}, 40, 185–220 (1984).

\bibitem{GK} K. Gawedzki, K., A. Kupiainen, Non-trivial continuum limit of a $\varphi_{4}^{4}$ model with negative coupling constant." Nuclear Physics B 257 (1985): 474-504.

\bibitem{GJ68} J. Glimm and A. Jaffe, A $\lambda \varphi_{2}^{4}$ quantum theory without cutoffs, I, Phys. Rev., 176, 1945–
1951 (1968), II: {\it Ann. Math.}, 91, 362–401 (1970), III: {\it Acta Math.}, 125, 203–267 (1970), IV: {\it J. Math. Phys.}, 13, 1568–1584 (1972).

\bibitem{GJ73} J. Glimm and A. Jaffe, Positivity of the $\varphi_{3}^{4}$ Hamiltonian, {\it Fort. Phys.}, 21, 327–376 (1973).

\bibitem{GJS74} J. Glimm, A. Jaffe and T. Spencer, The Wightman axioms and particle structure in the
$P(\varphi)_{2}$ quantum field model, {\it Ann. Math.}, 100, 585–632 (1974).

\bibitem{GJ74-2} J. Glimm and A. Jaffe, Remark on the existence $\varphi_{4}^{4}$, {\it Phys. Rev. Let.}, 33, 440–442 (1974).

\bibitem{GJ87} J. Glimm and A. Jaffe, {\it Quantum Physics: a Functional Integral Point of View}, Springer Verlag, 1987.

\bibitem{GRS75} F. Guerra, L. Rosen and B. Simon, The $P(\varphi)_{2}$ Euclidean quantum field theory as classical
statistical mechanics, {\it Ann. Math.}, 101, 111–259 (1975).

\bibitem{Haag} R. Haag, {\it Local Quantum Physics,} Springer Verlag 1992.

\bibitem{Hegerfeldt} G.C. Hegerfeldt, From Euclidean to relativistic fields and on the notion of Markoff
fields, {\it Commun. Math. Phys.}, 35, 155-171, 1974.

\bibitem{HK} R. Hoegh-Krohn, A General class of quantum fields without cut-offs in two space-time
dimensions, {\it Commun. Math. Phys.}, 21, 244–255 (1971).

\bibitem{Jaffe2007} A. Jaffe, Quantum field theory and relativity, in {\it Group representations, Ergodic Theory and Mathematical Physics: A tribute to
Georges W. Makey}, E.Doran, C.Moore and R. Zimmer (eds), Compt. Math, 449, pp. 209-246, 2008 (available at http://arthurjaffe.org).

\bibitem{Jost} R. Jost, {\it General Theory of Quantized Fields}, American Mathematical Society, Providence,
RI 1965.

\bibitem{Klauder} J.R. Klauder, {\it A Modern Approach to Functional Integration}, Birh\"auser 2011.

\bibitem{Kleinert-S} H. Kleinert, V. Schulte-Frohlinde, {\it Critical Properties of $\varphi^{4}$-Theories }, World Scientific, 2001.

\bibitem{LN} A. Lenard, C.M. Newman, Infinite volume asymptotics in $P(\varphi)_{2}$ field theory, {\it  Commun.
Math. Phys.}, 39, 243-250 (1974).

\bibitem{MS} J. Magnen and R. S\'en\'eor, The infinite volume limit of the $\varphi_{3}^{4}$ model, {\it  Ann. Inst. Henri
Poincar\'e}, 24, 95–159 (1976).

\bibitem{Malyshev} V.A. Malyshev, Probabilistic aspects of quantum field theory, Journ. Soviet. Math. 13 (4) 479-505, 1980.

\bibitem{Nelson1} E. Nelson, A quartic interaction in two dimensions, in {\it Mathematical Theory
of Elementary Particles }, R. Goodman and I. Segal, Eds., MIT Press, Cambridge, MA, 1966.

\bibitem{Nelson2} E. Nelson, The construction of quantum fields from Markov fields, {\it J. Funct. Analysis}, 12, 211-217 (1973).

\bibitem{Nelson3} E. Nelson, Probability theory and Euclidean field theory, in {\it Constructive Quantum Field
Theory}, G. Velo and A.S. Wightman Eds, Springer-Verlag, New York, pp. 94–124, 1973.

\bibitem{Neveu} J. Neveu, {\it Bases Mathématiques du Calcul des Probabilités}, 2d Ed. Masson.

\bibitem{Newman1} C. Newman, The construction of two-dimensionnal Markoff field with application to quantum field theory,
{\it J. Funct. Analysis} 14, 44-61 (1973).

\bibitem{Newman75} C.M. Newman, Inequalities for Ising models and field theories which obey the Lee-Yang theorem, 
{\it Commun. Math. Phys.} 41.1: 1-9 (1975)

\bibitem{OS} K. Osterwalder and R. Schrader, Axioms for Euclidean Green’s functions, {\it  Commun.
Math. Phys.}, 31, 83-112, 1973 and 42, 281–305 (1975).

\bibitem{Park} Y.M. Park, Convergence of lattice approximations and infinite volume limit in the  $(\lambda \varphi^{4}-\sigma \varphi^{2}-\mu \varphi)_{3} $ 
field theory, {\it J. Math. Phys.}, 18, 354–366 (1977).

\bibitem{PF} V.I.Pitebarg, V.R. Fatalov, The Laplace method for probability measures in Banach spaces, {\it Russ. Math. Surv.} Vo. 50, 6 1995

\bibitem{Rivasseau} V. Rivasseau, {\it From Perturbative to Constructive Renormalization} Princeton University
Press, Princeton, NJ) 1991.

\bibitem{Salmhofer} M. Salmhofer, {\it Renormalization: an introduction. Springer}, 2013

\bibitem{Schrader1} R. Schrader, On the Euclidean version of Haag's Theorem in $P(\varphi)_{2}$ theories, {\it  Commun.
Math. Phys.}, 36, 133-136 (1974). 

\bibitem{Schrader} R. Schrader, A possible constructive approach to $\varphi_{4}^{4}$, {\it  Commun.
Math. Phys.}, 49, 131-153 (1976). 

\bibitem{Simon} B. Simon, {\it The $P(\varphi)_{2}$ Quantum Field Theory}, Princeton Univ. Press, Princeton 1974.

\bibitem{SW} R.F. Streater and A.S. Wightman, {\it PCT, Spin and Statistics, and All That}, Reading,
Mass., Benjamin/Cummings Publ. Co. 1964.

\bibitem{Str} F. Strocchi, {\it Selected Topics on the General Properties of Quantum Field Theory}, World
Scientific, Singapore 1993.

\bibitem{Summers} S.J. Summers, A Perspective on Constructive Quantum FieldTheory, 2016;

\bibitem{Symanzik} K. Symanzik, Euclidean quantum field theory, I. Equations for a scalar model, {\it J. Math.
Phys.}, 7, 510–525, 1966.

\bibitem{Zinoviev} Y.M. Zinoviev, Equivalence of Euclidean and Wightman field theories, Comm. Math. Phys. 174, 1-27, 1995.

\end{thebibliography}
\end{document}